\documentclass[10pt]{article}
\usepackage{amsmath}
\usepackage{amssymb}
\usepackage{amsthm}
\usepackage{amscd}

\newtheorem{thm}{Theorem}[section]
\newtheorem{prop}[thm]{Proposition} 
\newtheorem{cor}[thm]{Corollary} 

\newtheorem{lem}[thm]{Lemma}
\newtheorem{defn}[thm]{Definition}

\theoremstyle{definition}
\newtheorem{rem}[thm]{Remark}

\newtheorem{exms}[thm]{Examples}

\begin{document}


\begin{center}
\textbf{\large{Reflexivity of non commutative Hardy Algebras}}
\end{center}

\begin{center}
Leonid Helmer
\end{center}

\begin{abstract}
Let $H^{\infty}(E)$ be a non commutative Hardy algebra, associated with a $W^*$-correspondence $E$. These algebras were introduced in 2004, ~\cite{MuS3}, by P. Muhly and B. Solel, and generalize the classical Hardy algebra of the unit disc $H^{\infty}(\mathbb{D})$.
As a special case one obtains also the algebra $\mathcal{F}^{\infty}$ of Popescu, which is $H^{\infty}(\mathbb{C}^n)$ in our setting.

In this paper we view the algebra $H^\infty(E)$ as acting on a Hilbert space via an induced representation. We write it $\rho(H^{\infty}(E))$ and we study the reflexivity of $\rho(H^{\infty}(E))$. This question was studied by A. Arias and G. Popescu in the context of the algebra $\mathcal{F}^{\infty}$, and by other authors in several other special cases. As it will be clear from our work, the extension to the case of a general $W^*$-correspondence $E$ over a general $W^*$-algebra $M$ requires new techniques and approach.

We obtain some partial results in the general case and we turn to the case of a correspondence over factor. Under some additional assumptions on the representation $\pi:M\rightarrow B(H)$ we show that $\rho_\pi(H^{\infty}(E))$ is reflexive. Then we apply these results to analytic crossed products $\rho(H^{\infty}(\ _{\alpha}M))$ and obtain their reflexivity for any automorphism $\alpha\in Aut(M)$ whenever $M$ is a factor. Finally, we show also the reflexivity of the compression of the Hardy algebra to a suitable coinvariant subspace $\mathfrak{M}$, which may be thought of as a generalized symmetric Fock space.
\end{abstract}

\section{Introduction.}

In this paper we consider the question of reflexivity of the non commutative Hardy algebras $\rho(H^\infty)(E)$.
These algebras, that were introduced by Muhly and Solel in \cite{MuS2} and \cite{MuS3}, can be viewed as a far reaching generalizations of the classical Hardy algebra $H^{\infty}(\mathbb{D})$.

Remember that $H^{\infty}(\mathbb{D})$ is the algebra of all the bounded analytic functions on the unit disc $\mathbb{D}\subset \mathbb{C}$, and may be identified with the WOT closed algebra, generated by the unilateral shift on the Hilbert space $l^2(\mathbb{Z}_+)$. In ~\cite{Po89} G. Popescu generalized it to WOT-closed algebras that are generated by $d$ shifts and denoted by $\mathcal{F}^{\infty}_d$. Note that the free semigroup algebra $\mathcal{L}_n$ studied in the late 90-th by K. Davidson and D. Pitts coincides with $\mathcal{F}^{\infty}$, (see $~\cite{DavP2}$).

These algebras were further generalized by Muhly and Solel who introduced in ~\cite{MuS2} and ~\cite{MuS3} the non commutative tensor algebras $\mathcal{T_+}(E)$ and Hardy algebras $H^{\infty}(E)$, associated with a $C^*$- or a $W^*$-correspondence $E$.
By a right Hilbert $C^*$-module over the $C^*$-algebra $A$ we mean a right $A$-module $E$, equipped with an $A$-valued inner product, that is a function $E\times E\rightarrow A$ such that $(\xi,\zeta)\mapsto\langle\xi,\zeta\rangle\in A$, for $\xi,\zeta\in A$, $a\in A$, and

1)$\langle\xi,\zeta a\rangle=\langle\xi,\zeta\rangle a$,

2)$\langle\xi,\zeta_1+\zeta_2 \rangle=\langle\xi,\zeta_1\rangle +\langle\xi,\zeta_2\rangle$,

3) $\langle\xi,\zeta \rangle=\langle\zeta,\xi\rangle^*$,

4) $\langle\xi,\xi \rangle\geq 0$, and $\langle\xi,\xi \rangle=0$ if and only if $\xi=0$.

This inner product defines a norm on $E$ by the formula $\|\xi\|^2:=\|\langle\xi,\xi\rangle\|$. The completion of $(E,\|\cdot\|)$ with respect to this norm is called a (right) Hilbert $C^*$-module.

By a bounded operator $T$ we mean a right module map $T:E\rightarrow E$, which is bounded in the norm defined above. Let $T$ be a map from Hilbert $C^*$-module $E$ into a Hilbert $C^*$-module $F$. We say that $T$ is
adjointable if there exists a map $T^*:F\rightarrow E$ such that

$$\langle Tx,y\rangle=\langle x,T^*y\rangle,$$
$x\in E$, $y\in F$. Such $T^*$ is called the adjoint map or the adjoint operator for $T$.
Every adjointable map is $A$-linear and bounded, (~\cite{Lan}).
The set of all adjointable maps from $E$ to $F$ we denote by $\mathcal{L}(E,F)=\mathcal{L}_A(E,F)$, and if $E=F$ we write $\mathcal{L}(E)$. It is known that $\mathcal{L}(E)$ is a $C^*$-algebra.
The standard reference to the theory of Hilbert $C^*$-modules is ~\cite{Lan}.

Every $C^*$-algebra $A$ has a natural structure of a Hilbert $C^*$-module over itself. The
right action is the multiplication in $A$ and for the inner product set $\langle a, b\rangle = a^*b$.
The Hilbert $C^*$-module over the $C^*$-algebra $A$ is called self-dual if every bounded $A$-linear map $T : E \rightarrow A$ is given by an inner product, which means that there is $x\in E$ such that for every
$y \in E$, $Ty = \langle x, y\rangle$. It is known that if $E$ is self-dual then $B(E) = \mathcal{L}(E)$ (~\cite{Lan}).

In this paper we are interested in the Hilbert modules over $W^*$-algebras.
Recall that by a $W^*$-algebra we mean a $C^*$-algebra $M$ that admits a faithful representations $\pi: M\rightarrow B(H)$ such that $\pi(M)$ is von Neumann algebra on the Hilbert space $H$. According to the abstract characterization given by Sakai, ~\cite{Sa1}, $C^*$-algebra $M$ is a $W^*$-algebra if and only if it is isomorphic to a dual Banach space.

A detailed study of  Hilbert modules over $W^*$-algebras was made by Paschke in ~\cite{Pa}.

A $C^*$-Hilbert module $E$ over the $W^*$-algebra $M$ which is also self-dual will be called a Hilbert $W^*$-module. In ~\cite{Pa} Paschke proved that every Hilbert $C^*$-module over the $W^*$-algebra $M$ admits the so called self-dual completion. In this case both $E$ and $\mathcal{L}(E)$ are dual spaces in the sense of Banach space theory, ~\cite{Pa}. In particular, $\mathcal{L}(E)$ is a $W^*$-algebra. The weak$^*$-topology on $E$, that comes from the structure of the dual space on $E$ will be called the $\sigma$-topology (following the terminology of ~\cite{BDH}).

By a $W^*$-correspondence we mean a right Hilbert $W^*$-module which is made into a bimodule over $M$ by some normal $*$-homomorphism $\phi$ of $M$ into the $W^*$-algebra of adjointable operators $\mathcal{L}(E)$.
Associated with a $W^*$-correspondence $E$ we have another $W^*$-correspondence $\mathcal{F}(E)$ over the same algebra $M$, that is defined to be the direct sum $M\oplus E\oplus E^{\otimes 2}\oplus...$ of the internal tensor powers of $E$. An exact definitions will be given in the next section. $\mathcal{F}(E)$ is called the full Fock correspondence and, in fact, it is a $W^*$-correspondence with the left action of $M$ denoted by $\phi_{\infty}$, which is a natural extension of $\phi$ to a representation of $M$ in the algebra of adjointable operators on $\mathcal{F}(E)$. The non commutative Hardy algebra of a correspondence $E$ is, by definition, the weak$^*$-closure in $\mathcal{L}(\mathcal{F}(E))$ of the algebra spanned by operators of the form $T_{\xi}$, $\xi\in E$, where $T_\xi(\eta):=\xi\otimes \eta$, $\eta\in \mathcal{F}(E)$, and $\phi_{\infty}(a)$, $a\in M$. In fact, Muhly and Solel defined this Hardy algebra as the weak$^*$ closure of the noncommutative tensor algebra $\mathcal{T}_+(E)$. The algebra $\mathcal{T}_+(E)$ was defined first in ~\cite{MuS2} as the norm closed (nonselfadjoint) algebra spanned by the same set of generators, and it generalizes the noncommutative disc algebra $\mathcal{A}_n$ of Popescu, which in turn is a noncommutative generalization of the classical disc algebra. Finally, the $C^*$-algebra, generated by the same set of operators is called the Toeplitz $C^*$-algebra of the given correspondence. In this context the algebra $\mathcal{F}^{\infty}$ coincides with the Hardy algebra $H^{\infty}(\mathbb{C}^n)$, where $\mathbb{C}^n$ is considered as the Hilbert $W^*$-module over $\mathbb{C}$.

These non commutative Hardy algebras, or simply Hardy algebras, generalize a wide class of known nonselfadjoint operator algebras. If we take$M=E=\mathbb{C}$ then $\mathcal{F}(E)$ is the Hilbert space $l^2(\mathbb{Z}_+)$, and the associated Hardy algebra is the classical algebra $H^{\infty}(\mathbb{D})$. Various choices of $W^*$-correspondence $E$ give us such algebras as the algebra $\mathcal{F}^{\infty}$ of G. Popescu, the free semigroup algebras, quiver algebras and analytic crossed products.

In this paper we view the algebra $H^{\infty}(E)$ as acting on a Hilbert space via an induced representation (induced in sense of Rieffel, see~\cite{Rie}). Given a representation $\pi$ of $M$ on $H$, the induced representation, written $\rho_\pi$ (or simply $\rho$,) is a representation of $H^{\infty}(E)$ on the Hilbert space $\mathcal{F}(E)\otimes_\pi H$ defined by sending $X\in H^\infty(E)$ to $X\otimes I_H$. Thus we consider the question of reflexivity of the algebra $\rho(H^{\infty}(E))$ and our results may be viewed as an extension of the results of A.Arias and G.Popescu in ~\cite{APo}.
It will be clear from our work that the extension to more general von Neumann algebra $M$ requires new techniques and approach. A key tool that we will need and use here is the concept of duality for $W^*$-correspondences, developed in ~\cite{MuS3}.

Let $\mathcal{A}$ be any algebra acting on a Hilbert space $H$. Then the  algebra $\mathcal{A}$ is said to be reflexive if it is defined by its invariant subspaces.
In a more details, with every subset $\mathcal{A}\subseteq B(H)$ let $Lat\ \mathcal{A}$ be the lattice of all $\mathcal{A}$-invariant closed subspaces in $H$. Thus,
$$Lat\ \mathcal{A}:=\{\mathcal{M}\subseteq H: A\mathcal{M}\subseteq
\mathcal{M}, \forall{A}\in \mathcal{A},\ \mathcal{M}\ \text{is a closed subspace in}\ H\}.$$

Let $\mathcal{L}$ be a set of closed subspaces in $H$. Then the operator $Alg$ associates to $\mathcal{L}$ an algebra $Alg\ \mathcal{L}\subseteq B(H)$ as follows.
$$Alg\ \mathcal{L}:=\{A\in B(H):\mathcal{L}\subseteq Lat\ A\},$$
i.e. $Alg\ \mathcal{L}=\{A\in B(H):A\mathcal{M}\subseteq \mathcal{M}, \forall{\mathcal{M}}\in \mathcal{L}\}$.
Clearly, $Alg\ \mathcal{L}$ is a unital WOT-closed (hence ultraweakly closed) subalgebra in $B(H)$.

Let $\mathcal{A}$ be a subalgebra in $B(H)$. Then always $\mathcal{A}\subseteq Alg\ Lat\mathcal{A}$, and by definition, $\mathcal{A}$ is reflexive if
$$\mathcal{A}=Alg\ Lat\ \mathcal{A}.$$
Thus, every reflexive algebra $\mathcal{A}$ is necessarily unital and ultraweakly closed.

A single operator $T$ is called reflexive if the WOT closed algebra, generated by $T$ and identity $I$, is reflexive. We denote this algebra by $W(T)$. The unilateral shift $S$ is reflexive. In this case $W(S)$ is the algebra of analytic Toeplitz operators, and the weak topology and the ultraweak topology coincide when restricted to $W(S)$. Thus the Hardy algebra
$H^{\infty}(\mathbb{T})\cong \rho(H^{\infty}(\mathbb{C}))$ is reflexive.

The following simple example shows that not every WOT closed algebra is reflexive. Let
$\mathcal{A}$ be an algebra of all $2\times 2$ matrices over $\mathbb{C}$ of the form
$\left(
    \begin{array}{cc}
      a & b \\
      0 & a \\
    \end{array}
  \right),\,\,\ a,b\in \mathbb{C}.
$
Clearly, $\mathcal{A}$ is WOT closed. But it is easy to see that
$$Alg\ Lat\ \mathcal{A}=\{\left(
                             \begin{array}{cc}
                               a & b \\
                               0 & c \\
                             \end{array}
                           \right):\, a,b,c\in \mathbb{C}
\}$$
Thus, $\mathcal{A}$ is not reflexive.
The notion of reflexivity was introduced first by H. Radjavi and P. Rosenthal in ~\cite{RadRos1} (the terminology was suggested by P. Halmos). It is easy to see that every von Neumann algebra is reflexive, which is equivalent to the von Neumann bicommutant theorem.
Since the paper ~\cite{RadRos1} appears, the subject of reflexivity was a subject of intensive study and generalizations. A good general overview of reflexivity is given by D. Hadwin in ~\cite{Hadw}.

The works which are closely connected to our theme are ~\cite{APo}, ~\cite{Po6}, ~\cite{Dav1}, ~\cite{Dav2}, ~\cite{DavP2}, \cite{Ken}, ~\cite{KriPow} ~\cite{Ka}. In particular, in ~\cite{APo}, A. Arias and G. Popescu proved reflexivity of the algebra $\mathcal{F}^{\infty}$, which, as we shall see later, coincides with $\rho(H^{\infty}(\mathbb{C}^{n})$ in our setting. Later, in ~\cite{Po6} Popescu proved reflexivity of this algebra when it is compressed to the symmetric Fock space. This algebra also was studied by Davidson and Pitts and by Davidson in the context of the free semigroup algebra $\mathcal{L}_n$, see ~\cite{Dav1}, ~\cite{Dav2}, ~\cite{DavP2}. In the recent work `\cite{Ken}, M. Kennedy showed the reflexivity of all the free semigroup algebras and the hyperreflexivity of some of them (the definition of
hyperreflexivity will be given later). In the work ~\cite{KriPow}, D. Kribs and S. Power initiated the study of the free semigroupoid algebras $\mathcal{L}_G$ and in particular proved their reflexivity. In ~\cite{Ka} E. Kakariadis showed the reflexivity of the analytic crossed product $\rho(H^\infty( _{\alpha}M))$ in the special case when the $\alpha$ is a unitary implemented automorphism of von Neumann algebra.

This paper is organized as follows. After introducing preliminaries we start with some simple observations about hyperreflexivity and then turn to our main results. One of the our central result is Theorem $\ref{Repres. of Phi_j}$ on the matrix representation of Fourier coefficients. To show reflexivity of $\rho(H^{\infty}(E))$ we need to show that for every $Y\otimes I_H\in Alg\ Lat\ \rho(H^{\infty}(E))$ each Fourier coefficient $\{\Phi_j(Y)\otimes I_H\}$, $j=0,1,2...$ is in $\rho(H^{\infty}(E))$. This approach generalizes the approach of Popescu and Arias in their proof of the reflexivity of $\mathcal{F}^{\infty}$, ~\cite{APo}.
To this end, we introduce the coinvariant subspace $\mathfrak{M}$ that may be thought of as a generalization of the symmetric Fock space and seems to be of particular interest. Using this subspace we are able to show that $\Phi_0(Y)$ is in $\rho(H^{\infty}(E))$. For the Fourier coefficients $\Phi_j(Y)$ with $j\geq 1$ of an arbitrary $Y\in Alg\ Lat\ (H^{\infty}(E))$, this method does not work for a general $W^*$-algebra. But if $M$ is assumed to be a factor, then under an additional assumption on the representation $\pi$ we are able to show the reflexivity of our algebra. As an example we consider the analytic crossed products $\rho(H^{\infty}( _{\alpha}M))$. Finally, we consider the compression of the Hardy algebra to the subspace $\mathfrak{M}$, i.e. we consider the algebra $Q\rho(H^{\infty}(E))|_{\mathfrak{M}}$, where $Q$ is the projection onto $\mathfrak{M}$, and show that it is reflexive.

This work is based on part of the author’s Ph.D. thesis.

\section{Preliminaries.}
\subsection{$W^*$-correspondences and Hardy algebras.}

We start by recalling the definition of a $W^*$-correspondence.
\begin{defn} Let $E$ be a (right) Hilbert $W^*$-module over a $W^*$-algebra $M$, that is a self-dual $C^*$-module over the $W^*$-algebra $M$, and let $\phi: N\rightarrow \mathcal{L}(E)$ be a normal $*$-homomorphism of the $W^*$-algebra $N$ into the algebra of adjointable operators $\mathcal{L}(E)$. Then $\phi$ defines on $E$ the structure of the left module over $N$. This $N$-$M$-bimodule is called the $W^*$-correspondence from $N$ to $M$. If $N=M$ we speak about a $W^*$-correspondence over $M$.
\end{defn}

In what follows we always assume that our $W^*$-correspondence over $M$ is essential as a left $M$-module,
meaning that $\phi(M)E$ is dense in $E$ in the $\sigma$-topology.

The obvious example of $W^*$-correspondence over $M=\mathbb{C}$ is the ordinary Hilbert space $H$ with the inner product is taken to be linear in the second variable.

Let $M$ be a $W^*$-algebra which we view as the Hilbert $C^*$-module over itself. The self dual completion gives rise to a $W^*$-module. If $\alpha$ is some normal automorphism of $M$ we set $\phi(a)b:=\alpha(a)b$ for the left action. Then $M$ turns out to be a $W^*$-correspondence over itself and is denoted by $_{\alpha}M$. More generally, if $\alpha: M\rightarrow M$ is a normal $*$-homomorphism, that is
$\alpha\in End(M)$, then $E:=\overline{\alpha(M)M}$ turns out to be a $M^*$-correspondence with $\alpha$ as the left action. Note that $\alpha(1)=p$ - some projection of $M$. Thus, we conclude that $E$ is a cyclic right module of the form $pM$. Later, we will return to this module in more details.

Let $M$, $N$ and $Q$ be three $W^*$-algebras, $E$ be a $M$-$N$ $W^*$- correspondence, and $F$ be a $N$-$Q$ $W^*$-correspondence. Write $\pi$ for the left action of $N$ on $F$, and write $E\otimes_{alg} F$ for the algebraic tensor product of $\mathbb{C}$-vector spaces
$E$ and $F$. By the internal $C^*$-tensor product (balanced over $\pi$) of these modules we mean the  Hausdorff completion of $E\otimes_{alg} F$ by the inner product defined by
$\langle
\xi_1\otimes \eta_1,\xi_2\otimes
\eta_2\rangle=\langle\eta_1,\pi(\langle\xi_1,\xi_2\rangle)\eta_2\rangle$. This tensor product will be denoted by $E\otimes_\pi F$. For the right action of $Q$ we set $(\xi\otimes\eta)
b=\xi\otimes(\eta b)$, $b\in Q$. The self-dual completion of $E\otimes_{\pi}F$ gives rise to a Hilbert $W^*$-module over the $W^*$-algebra $Q$. For the left action we set
$\phi_{E\otimes_\pi F}(a)(\xi\otimes \eta)=(\phi_E(a)\xi)\otimes \eta$. Thus, we obtain on $E\otimes_{\pi}F$ a structure of a $W^*$-correspondence from $M$ to $Q$.

It is easy to see that if $E$ is a Hilbert module over a $W^*$-algebra $M$ and $\pi:M\rightarrow B(H)$  a normal representation of $M$ on the Hilbert space $H$ then $E\otimes_{\pi}H$ is Hilbert space.

Let $E$ be a $W^*$-correspondence over a $W^*$-algebra $M$ with a left action defined as usual by a normal $*$-homomorphism $\phi$.
For each $n\geq 0$ let $E^{\otimes n}$ be the self-dual internal tensor power (balanced over $\phi$). So, $E^{\otimes n}$ itself turns out to be a $W^*$-correspondence in a natural way, with the left action $\xi \mapsto \phi_n(a)\xi=(\phi(a)\xi _{1})\otimes ...\otimes \xi_n$, $\xi=\xi_1\otimes...\otimes\xi_n\in E^{\otimes n}$, and with an $M$-valued inner product which comes from the construction of the internal tensor product.

Let $\{E_i:i\in I\}$ be a family of $W^*$-correspondences over $M$. By $\sum_i^{\oplus}E_i$ we denote
the set of all sequences $x=(x_i)$, $x_i\in E_i$, such that $\sum_i\langle x_i,x_i\rangle$
converges in $M$ (considered as a $C^*$-algebra). For $x=(x_i)$ and $y=(y_i)$ in $\sum_i^{\oplus}E_i$, we define $\langle
x,y\rangle =\sum_i\langle x_i,y_i\rangle$.
This defines an inner product on $\sum_i^{\oplus}E_i$ and in fact $\sum_i^{\oplus}E_i$ is complete in
the norm defined by this inner product. With the obvious right action of $M$ on $\sum_i^{\oplus}E_i$, this module is a Hilbert module over $M$ in the $C^*$-sense. Then the self-dual completion of $\sum_i^{\oplus}E_i$ will be called the ultraweak direct sum of the family $\{E_i:i\in I\}$ (for the explicit description see ~\cite{Pa}. For the left action we set $\phi=\sum_i^{\oplus}\phi_{E_i}$, where $\phi_{E_i}$ is a left action of $M$ on $E_i$. Clearly, we obtain on $\sum_i^{\oplus}E_i$ a structure of the $W^*$-correspondence over $M$.

We form the full Fock space $\mathcal{F}(E)=\sum^{\oplus}_{n\geq 0}E^{\otimes n}$, where $E^{\otimes 0}=M$ and the direct sum taken in the ultraweak sense. This is a $W^*$-correspondence with left action given by
$\phi _{\infty}:M\rightarrow \mathcal{L}(\mathcal{F}(E))$, where $\phi _{\infty}(a)=\sum _{n\geq 0} \phi_n(a)$. The $M$-valued inner product on $\mathcal{F}(E)$ is defined in an obvious way.

For each $\xi\in E$ and each $\eta\in \mathcal{F}(E)$, let $T_{\xi}: \eta \mapsto \xi \otimes \eta$ be a creation operator on
$\mathcal{F}(E)$. Clearly, $T_{\xi}\in \mathcal{L}(\mathcal{F}(E))$.

\begin{defn}\label{Defn of tensor and Hardy alg} Given a $W^*$-correspondence $E$ over a $W^*$-algebra $M$.

1) The norm closed subalgebra of $\mathcal{L}(\mathcal{F}(E))$, generated by all creation operators $T_{\xi}$, $\xi\in E$, and all operators $\phi _{\infty}(a)$,  $a\in M$, is called the tensor algebra of $E$. It is denoted by $\mathcal{T}_{+}(E)$.

2) The Hardy algebra $H^{\infty}(E)$ is the ultraweak closure
of $\mathcal{T}_{+}(E)$ in the $W^*$-algebra $\mathcal{L}(\mathcal{F}(E))$.
\end{defn}

\begin{exms}\label{Examples of tensor and Hardy alg,classical and Popescu}
(1) Let $A=E=\mathbb{C}$. Then $\mathcal{F}(E)=l^2(\mathbb{Z}_+)$ and the algebra
$\mathcal{T}_{+}(E)$ is the algebra of analytic Toeplitz operators with continuous symbols, so it can be identified with the disc algebra $A(\mathbb{D})$. The algebra
$H^{\infty}(E)$, in this case, is $H^{\infty}(\mathbb{D})$.

(2) Let $A=\mathbb{C}$ and let $E$ be an $n$-dimensional Hilbert space over $\mathbb{C}$, i.e., $H=\mathbb{C}^n$.
In this case $\mathcal{T}_+(E)$ is the non commutative disc algebra $\mathcal{A}_n$, studied by Popescu and others. The algebra $H^{\infty}(\mathbb{C}^n)$ was denoted $\mathcal{F}^{\infty}_n$ by Popescu. This algebra can be identified with the free semigroup algebra $\mathcal{L}_n$ studied by Davidson and Pitts.
\end{exms}

Let $\pi:M\rightarrow B(H)$ be a normal representation of a $W^*$-algebra $M$ on a Hilbert space $H$ and let $E$ be a $W^*$-correspondence over $M$. As we already noted, the $W^*$-internal tensor product $E\otimes_{\pi}H$ is a Hilbert space. The representation $\pi^{E}:\mathcal{L}(E)\rightarrow B(E\otimes _{\pi}H)$ defined by

$$\pi^{E}: S\mapsto S\otimes I_H, \,\,\,\ \forall S\in \mathcal{L}(E).$$
is called the induced representation (in the sense of Rieffel). If $\pi$ is a faithful normal representation then $\pi^E$ maps $\mathcal{L}(E)$ into
$B(E\otimes _{\pi}H)$ homeomorphically with respect to the ultraweak topologies, ~\cite[Lemma 2.1]{MuS3}.

In this work, we consider the image of $H^{\infty}(E)$ under an induced representation, defined as follows. Let $\pi:M\rightarrow B(H)$ be a faithful normal representation. For a $W^*$-correspondence $E$ over $M$ let
$\pi^{\mathcal{F}(E)}$ be the induced representation of $\mathcal{L}(\mathcal{F}(E))$ in $B(\mathcal{F}(E)\otimes_{\pi}H)$.
Then the induced representation of the Hardy algebra $H^{\infty}(E)$ is the restriction
\begin{equation}\label{Induced repr of Hardy alg}
\rho:=\pi^{\mathcal{F}(E)}|_{H^{\infty}(E)}:H^{\infty}(E)\rightarrow B(\mathcal{F}(E)\otimes_{\pi}H).
\end{equation}
This restriction is an ultraweakly continuous representation
of ${H^{\infty}(E)}$ and the image $\rho(H^{\infty}(E))$ is an ultraweakly closed subalgebra of $B(\mathcal{F}(E)\otimes_{\pi}H)$. We shall refer to $\rho$ as the representation induced by $\pi$. Later, when we discuss several representation of $H^{\infty}(E)$ that are induced by different representations $\pi$, $\sigma$ etc. of $M$, we shall write $\rho_{\pi}$, $\rho_\sigma$ etc.

So, $\rho(H^{\infty}(E))$ acts on
$\mathcal{F}(E)\otimes_{\pi}H$ and $\rho$ is defined by
$$\rho:X\mapsto X\otimes I_H,\,\,\,\forall{X}\in H^{\infty}(E).$$

Note that the notion of the induced representation generalizes the notion of pure isometry in the theory of a single operator.

We will frequently use the following result of Rieffel ~\cite[Theorem 6.23]{Rie}.
The formulation here is in a form convenient for us (~\cite[p. 853]{MuS1}).
\begin{thm}\label{Rieffel thm}. Let $E$ be a $W^*$-correspondence over the algebra $M$ and $\pi:M\rightarrow B(H)$ be a normal faithful representation of $M$ on the Hilbert space $H$. Then the operator $R$ in $B(E\otimes_{\pi}H)$ commutes with $\pi^{E}(\mathcal{L}(E))$ if and only if $R$ is of the form $I_E\otimes X$, where $X\in \pi(M)'$, i.e., $\pi^{E}(\mathcal{L}(E))'=I_E\otimes \pi(M)'$.
\end{thm}

\subsection{Covariant representations.}

\begin{defn}
Let $E$ be a $W^*$-correspondence over a $W^*$-algebra $M$.

(1) By a covariant representation of $E$, or of the pair $(E,M)$, on a Hilbert space $H$, we mean
a pair $(T, \sigma)$, where $\sigma:M\rightarrow B(H)$ is a nondegenerate normal $*$-homomorphism, and
$T$ is a bimodule (with respect to $\sigma$) map $T:E\rightarrow B(H)$, that is a linear map such
that $T(\xi a)=T(\xi)\sigma(a)$ and $T(\phi(a)\xi)=\sigma(a)T(\xi)$, $\xi\in E$ and $a\in M$. We require also that $T$ will be continuous with respect to the $\sigma$-topology on $E$ and the ultraweak topology on $B(H)$.

(2) The representation $(T,\sigma)$ is called (completely) bounded, (completely) contractive, if so is the map $T$. For a completely contractive covariant representation we write also c.c.c.r.

(3) The covariant representation $(T,\sigma)$ is called isometric covariant representation
(i.c.r.) if $T(\xi)^*T(\eta)=\sigma(\langle \xi,\eta\rangle)$.

\end{defn}

The operator space structure on $E$ to which this definition refers is the one which comes from the embedding of $E$ into its so-called linking algebra $\mathfrak{L}(E)$, see ~\cite{MuS2}.

Every isometric covariant representation $(V,\pi)$ of $E$ is completely contractive (see ~\cite[ Corollary 2.13]{MuS2}).

As an important example let $\rho=\pi^{\mathcal{F}(E)}|_{H^{\infty}(E)}$ be an induced representation of the Hardy algebra $H^{\infty}(E)$. For the representation $\sigma$ set
$$\sigma=\pi^{\mathcal{F}(E)}\circ\phi_{\infty},$$
and set
$$V(\xi)=\pi^{\mathcal{F}(E)}(T_{\xi}),\,\,\,\ \xi\in E.$$
\begin{defn}\label{Def of ind cov rep}
The pair $(V,\sigma)$ is called the covariant representation induced by $\pi$, or simply the induced covariant representation (associated with $\rho$).
\end{defn}

It is easy to check that $(V,\sigma)$ in the above Definition is isometric, hence, is completely contractive.

Let $(T, \sigma)$ be a c.c.c.r. of $(E,M)$ on the Hilbert space $H$ as above.
With each such representation we associate the operator $\tilde{T}:E\otimes _{\sigma}H\rightarrow H$, that on the elementary tensors is defined by
$$\tilde{T}(\xi \otimes h):=T(\xi)(h).$$
$\tilde{T}$ is well defined since $T(\xi a)=T(\xi)\sigma(a)$.
In ~\cite{MuS2} Muhly and Solel show that the properties of $\tilde{T}$ reflect the properties of the covariant representation $(T,\sigma)$.
They proved that $(\alpha)$ $ \tilde{T}$ is bounded iff $T$ is
completely bounded, and in this case $\|T\|_{cb}=\|\tilde{T}\|$; $(\beta)$ $\tilde{T}$ is
contractive iff $T$ is completely contractive; and $(\gamma)$ $\tilde{T}$ is an isometry  iff
$(T,\sigma)$ is an isometric representation.
A simple calculation gives us the intertwining relation
\begin{equation}\label{intetwining for T}
\tilde{T}\sigma^{E}\circ\phi(a)=\tilde{T}(\phi(a)\otimes I_H)=\sigma (a)\tilde{T}, \,\,\ \forall
a\in A.
\end{equation}

We have the following lemma, taken from ~\cite[Lemma 2.5]{MuS3}.

\begin{lem}
There exists a bijective correspondence $(T,\sigma)\leftrightarrow \tilde{T}$, between all
completely contractive representations $(T,\sigma)$ of $E$ on a Hilbert space $H$, and all
contractive operators $\tilde{T}:E\otimes_{\sigma}H\rightarrow H$, that satisfy the relation
$\tilde{T}\sigma^{E}\circ\phi(a)=\sigma (a)\tilde{T}$, $\forall a\in A$.
Given a contractive operator $\tilde{T}:E\otimes_{\sigma}H\rightarrow H$ that satisfies the above
intertwining relation, then the associated covariant representation $(T,\sigma)$, is defined by
$T(\xi)h:=\tilde{T}(\xi\otimes h)$, $h\in H$ and $\xi\in E$.
\end{lem}

\begin{rem} Let $E$ be a $W^*$-correspondence over the algebra $M$ and let
$(T,\sigma)$ be a c.c.c.r. of $(E,M)$ on a Hilbert space $H$. It is shown in ~\cite{MuS3} that for every such c.c.c.r. there exists a completely contractive representation $\rho:\mathcal{T}_+(E)\rightarrow B(H)$ such that
$\rho(T_{\xi})=T(\xi)$ for every $\xi\in E$ and $\rho(\phi_{\infty}(a))=\sigma(a)$ for every $a\in M$. Moreover, the correspondence $(T,\sigma)\leftrightarrow\rho$ is a bijection between the set of all c.c.c.r. of $E$ and all completely contractive representations of $\mathcal{T}_+(E)$ whose restrictions to $\phi_{\infty}(M)$ are continuous with respect to the ultraweak topology on $\mathcal{L}(\mathcal{F}(E))$.

The representation $\rho$ of $\mathcal{T}_+(E)$ that corresponds to the c.c.c.r. $(T,\sigma)$ is called the integrated form of $(T,\sigma)$ and denoted by $\sigma\times T$. In its turn, the c.c.c.r. $(T,\sigma)$ is called the desintegrated form of $\rho$.
Preceding results show that, given a normal representation $\sigma$ of
$M$, the set of all completely contractive representations of $\mathcal{T}_+(E)$ whose restrictions to $\phi_{\infty}(M)$ is given by $\sigma$ can be
parameterized by the contractions $\tilde{T}\in B(E\otimes_{\sigma}H, H)$, that satisfy the relation $(\ref{intetwining for T})$.

In ~\cite{MuS3} it was shown that, if
the c.c.c.r. $(T,\sigma)$ of $(E,M)$ is such that
$\|\tilde{T}\|<1$, then the integrated form $\sigma\times T$
extends from $\mathcal{T}_+(E)$ to an ultraweakly continuous representation of $H^{\infty}(E)$. For a general c.c.c.r. $(T,\sigma)$, the question when such extention is possible $\sigma\times T$ is more delicate, see about this ~\cite{MuS8}.

In the above notations, the induced representation $\rho_\pi$ is an integrated form of the $(V,\sigma)$, the covariant induced representation of $E$ from Definition $\ref{Def of ind cov rep}$.
\end{rem}

We shall use the following notation. Let $(V,\sigma)$ be an isometric covariant representation of a general $W^*$-correspondence $E$ on a Hilbert space $G$. For every $n\geq 1$ write $(V^{\otimes n},\sigma)$ for the isometric covariant representation of $E^{\otimes n}$ on the same space $G$ defined by the formula $V^{\otimes n}(\xi_1\otimes...\otimes\xi_n)=V(\xi_1)\cdots V(\xi_n)$, $n\geq 1$. The associated isometric operator $\tilde{V}_n:E^{\otimes n}\otimes_\sigma G\rightarrow G$ (which is called the generalized power of $\tilde{V}$), satisfies the identity
$\tilde{V}_n\sigma^{E^{\otimes n}}\circ \phi_n= \tilde{V}_n (\phi_n\otimes I_{G_0})=\sigma \tilde{V}_n$. In this notation $\tilde{V}=\tilde{V}_1$.

With $(V,\sigma)$ we may associate the "shift" $\mathfrak{L}$, that acts on the lattice of $\sigma (M)$-invariant subspaces of $G$, and is defined as follows. Let $\mathcal{M}\in Lat(\sigma (M))$, then we set

\begin{equation}\label{DefnOfGeneralizedShift}
\mathfrak{L}(\mathcal{M}):=\bigvee \{V(\xi)k : \xi \in E, k\in \mathcal{M}\}.
\end{equation}
The $s$-power $\mathfrak{L}^s(\mathcal{M})$ is defined in the obvious way (with $\mathfrak{L}^0(\mathcal{M})=\mathcal{M})$).

The subspace $\mathcal{M}\in Lat(\sigma(M))$, as well as its projection $P_{\mathcal{M}}\in \sigma(M)'$, is called wandering with respect to $(V,\sigma)$, if the subspaces  $\mathfrak{L}^s(\mathcal{M})$, $s=0,1,...$, are mutually orthogonal.
Write $\sigma'$ for the restriction $\sigma|_{\mathcal{M}}$, where $\mathcal{M}$ is wandering. Then the Hilbert space $E^{\otimes s}\otimes_{\sigma'}\mathcal{M}$ is isometrically isomorphic (under the generalised power $\tilde{V}_s$) to $\mathfrak{L}^s(\mathcal{M})$. Hence, we obtain an isometric isomorphism
$$\mathcal{F}(E)\otimes_{\sigma'}\mathcal{M}\cong \sum ^{\oplus}_{s\geq 0}\mathfrak{L}^s(\mathcal{M}).$$

\subsection{Duality of $W^*$-correspondences and commutant.}
The principal tool that will be used often in this work is the duality of $W^*$-correspondences, that was developed in ~\cite[Section 3]{MuS3}. We shall need the notion of isomorphic $W^*$-correspondences. Let $E$ and $F$ be $W^*$-correspondences over $W^*$-algebras $M$ and $N$ respectively. The left action of $M$ on $E$ will be denoted as usual by $\phi$ and the left action of $N$ on $F$ by $\psi$, thus, $\psi:N\rightarrow \mathcal{L}(F)$ is a normal $^*$-homomorphism.
\begin{defn}\label{Isomorphism of W^*-corresp} An isomorphism of $E$ and $F$ is a pair $(\sigma,\Phi)$ where

1) $\sigma:M\rightarrow N$ is an isomorphism of $W^*$-algebras and

2) $\Phi:E\rightarrow F$ is a vector space isomorphism preserving the $\sigma$-topology, and which is also (a) a bimodule map, $\Phi(\phi(a)xb)=\psi(\sigma(a))\Phi(x)\sigma(b)$, $x\in E$ $a,b\in M$, and

(b) $\Phi$ "preserves" the inner product, $\langle\Phi(x),\Phi(y)\rangle=\sigma(\langle x,y\rangle)$, $x,y\in E$.
\end{defn}

Let $\pi: M\rightarrow B(H)$ be a normal representation of $M$ on a Hilbert space $H$. We put
\begin{equation}
E^{\pi}:=\{\eta:H\rightarrow E\otimes_{\pi}H:\eta \pi(a)=(\phi(a)\otimes
I_H)\eta, a\in M\}.
\end{equation}

On the set $E^{\pi}$ we define the structure of a $W^*$-correspondence over the von Neumann algebra $\pi(M)'$ putting
$\langle \eta,\zeta \rangle:=\eta^*\zeta$ for the $\pi(M)'$-valued inner product, $\eta,\zeta\in E^{\pi}$. It is easy to check that $\langle \eta,\zeta \rangle\in \pi(M)'$. For the bimodule operations: $X\cdot \eta=(I\otimes X)\eta$, and
$\eta\cdot Y=\eta Y$, where $X,Y \in \pi(M)'$.
\begin{defn} The $W^*$-correspondence $E^{\pi}$ is called the $\pi$-dual of $E$.
\end{defn}

Write $(E^{\pi})^*$ for the space of adjoints of the operators in $E^{\pi}$. By
$\overline{\mathbb{D}(E^{\pi})}$ and $\mathbb{D}(E^{\pi})$ we denote the (norm) closed and
open unit balls in $E^{\pi}$ respectively. Let $\eta^*\in
\overline{\mathbb{D}(E^{\pi})^*}$ and let $(T,\pi)$ the associated c.c.c.r. of $(E,M)$ on $H$ such that $\eta^*=\tilde{T}$. Since every $(T,\pi)$ satisfies relation
$(\ref{intetwining for T})$, we obtain that all the representations $(T,\pi)$ of $(E,M)$ are parameterized
by the points of $\overline{\mathbb{D}(E^{\pi})^*}$. Hence, all the completely contractive representations $\rho$
of $\mathcal{T}_+(E)$ such that $\rho \circ \phi_{\infty}=\pi$ are parameterized bijectively by
$\overline{\mathbb{D}(E^{\pi})^*}$.

Let $\iota:\pi(M)'\rightarrow B(H)$ be the identity representation.
Then we can form $E^{\pi,\iota}:=(E^{\pi})^{\iota}$.
So, $E^{\pi,\iota}=\{S:H\rightarrow E^{\pi}\otimes_{\iota}H:
S\iota(a)=\iota^{E^{\pi}}\circ\phi_{E^{\pi}}(a)S,
a\in \pi(M)'\}$. This is a $W^*$-correspondence over $\pi(M)''=\pi(M)$.

In ~\cite{MuS3} it was proved that for every faithful normal representation $\pi$ of a $W^*$-algebra $M$, every $W^*$-correspondence $E$ over $M$ is isomorphic to
$E^{\pi,\iota}$.
We give a short description of this isomorphism.

Let $L_{\xi}:h\mapsto \xi\otimes h$, $h\in H$,$\xi \in E$.
$L_{\xi}$ is a bounded linear map since $\|L_{\xi}h\|^2\leq \|\xi\|^2\|h\|^2$
and $L_{\xi}^*(\zeta\otimes h)=\pi(\langle \xi,\zeta \rangle)h$.
For each $\xi\in E$ we define the map $\hat{\xi}:H\rightarrow E^{\pi}\otimes_{\iota}H$ by
means of its adjoint:
$$\hat{\xi}^*(\eta\otimes h)=L_{\xi}^*(\eta(h)),$$
$\eta\otimes h\in E^{\pi}\otimes_{\iota} H$.

\begin{thm}(~\cite[Theorem 3.6]{MuS3})\label{Isomorphism with the second dual corresp}
If the representation $\pi$ of $M$ on $H$ is faithful, then the map $\xi\mapsto \hat{\xi}$ just defined, is an isomorphism of the $W^*$-correspondences $E$ and $E^{\pi,\iota}$.
\end{thm}

For every $k\geq 0$, let $U_k:E^{\otimes k}\otimes_{\pi}H\rightarrow (E^{\pi})^{\otimes k}\otimes_{\iota}H$
be the map defined in terms of its adjoint by
$U_k^*(\eta_1\otimes...\otimes \eta_n\otimes h)=(I_{E^{\otimes k-1}}\otimes \eta_1)...(I_{E}\otimes
\eta_{k-1})\eta_k(h)$.
It is proved in ~\cite{MuS3} that $U_k$
is a Hilbert space isomorphism from $E^{\otimes k}\otimes_{\pi}H$ onto $(E^{\pi})^{\otimes k}\otimes_{\iota}H$.

By Theorem $\ref{Isomorphism with the second dual corresp}$, for every $k\geq 1$ the $W^*$-correspondence $E^{\otimes k}$ over $M$ is isomorphic to the $W^*$-correspondence $(E^{\otimes k})^{\pi,\iota}\cong (E^{\pi,\iota})^{\otimes k}$. If $\xi\in E^{\otimes k}$ then the corresponding element $\widehat{\xi}\in (E^{\otimes k})^{\pi,\iota}$
is defined now by the formula
$$\widehat{\xi}^*(\eta_1\otimes...\otimes\eta_k\otimes h)=L_{\xi}^*U^*_k(\eta_1\otimes...\otimes\eta_k\otimes h),$$
where $L_{\xi}:h\mapsto \xi\otimes h$ is a bounded linear map from $H$ to $E^{\otimes k}\otimes_{\pi}H$.
Thus, we obtain
\begin{equation}\label{Formula for hat(xi)}
\hat{\xi}=U_kL_{\xi},\,\,\ \text{for}\, \xi\in E^{\otimes k}.
\end{equation}

For the dual correspondence ($\pi$-dual to $E$) we can form the (dual) Fock space $\mathcal{F}(E^{\pi})$, which is a $W^*$-correspondence over $\pi(M)'$, and the Hilbert space $\mathcal{F}(E^{\pi})\otimes_\iota H$. Let us define
$U:=\sum^{\oplus}_{k\geq 0}U_k$. It follows that
the map $U:=\sum^{\oplus}_{k\geq 0}U_k$ is a Hilbert space isomorphism from
$\mathcal{F}(E)\otimes_{\pi}H $ onto $\mathcal{F}(E^{\pi})\otimes_{\iota}H$, and its adjoint acts on decomposable tensors by
$U^*(\eta_1\otimes...\otimes\eta_n\otimes h)=(I_{E^{\otimes n-1}}\otimes \eta_1)...(I_{E}\otimes
\eta_{n-1})\eta_nh$.

\begin{defn}\label{Fourier trans definition}
The map $U_{\pi}=U:\mathcal{F}(E)\otimes_{\pi}H\rightarrow\mathcal{F}(E^{\pi})\otimes_{\iota}H$
will be called the Fourier transform determined by $\pi$.
\end{defn}

Let $\pi:M\rightarrow B(H)$ be a faithful normal representation. Then there exists a natural
canonical isometric representation of $E^{\pi}$ on $\mathcal{F}(E)\otimes_{\pi}H$ induced
by  $\pi$. Let $\nu:\pi(M)'\rightarrow B(\mathcal{F}(E)\otimes_{\pi}H)$ be a $*$-representation defined
by $\nu(b)=I_{\mathcal{F}(E)}\otimes b$. Then $\nu$ is a faithful normal
representation of the von Neumann algebra $\pi(M)'$. By Rieffel's Theorem $\ref{Rieffel thm}$,
$\pi^{\mathcal{F}(E)}(\mathcal{L}(\mathcal{F}(E)))'=\nu(\pi(M)')=\{I_{\mathcal{F}(E)}\otimes b :b\in
\pi(M)'\}$. Let $\eta\in E^{\pi}$. Then for every $n\geq 0$, the operators
$L_{\eta,n}:E^{\otimes n}\otimes_{\pi}H\rightarrow E^{\otimes n+1}\otimes_{\pi}H$ are defined
by $L_{\eta,n}(\xi\otimes h)=\xi\otimes \eta h$, where we have identified $E^{\otimes
n+1}\otimes_{\pi}H$ with $E^{\otimes n}\otimes_{\pi^{E}\circ \phi}(E \otimes_{\pi}H)$.
Since $\|L_{\eta,n}\|\leq \|\eta\|$, we may define the operator
$\Psi(\eta):\mathcal{F}(E)\otimes_{\pi}H\rightarrow \mathcal{F}(E)\otimes_{\pi}H$
by $\Psi(\eta)=\sum^{\oplus}_{k\geq 0}L_{\eta,k}$. Thus we may think of $\Psi(\eta)$ as $I_{\mathcal{F}(E)}\otimes \eta$ on
$\mathcal{F}(E)\otimes_{\pi}H$. It is easy to see that $\Psi$ is a bimodule map, and not hard to check that $(\Psi,\nu)$ is an isometric covariant representation of $(E^{\pi},\pi(M)')$ on the Hilbert space $\mathcal{F}(E)\otimes_{\pi}H$, (for more  details see ~\cite{MuS3}).

Now, combining the integrated form $\nu \times \Psi$ of $(\Psi,\nu)$ with the definition of the Fourier transform $U=U_\pi$
we obtain the formulas

\begin{equation}\label{Psi(eta)formula}
U^*\iota^{\mathcal{F}(E^{\pi})}(T_{\eta})U=\Psi(\eta),
\end{equation}
where $\eta\in E^{\pi}$ and $T_{\eta}$ is the corresponding creation operator
in $H^{\infty}(E^{\pi})$, and
\begin{equation}\label{pi_Psi(a)formula}
U^*\iota^{\mathcal{F}(E^{\pi})}(\phi_{E^{\pi},\infty}(b))U=\nu(b),
\end{equation}
where $b\in \pi(M)'$ and $\phi_{E^{\pi},\infty}$ is the left action of $\pi(M)'$ on $\mathcal{F}(E^{\pi})$.  This equality can be rewritten as
\begin{equation}\label{Intertwining M' and U}
U(I_{\mathcal{F}(E)}\otimes b)=(\phi_{E^{\pi},\infty}(b)\otimes I_H)U.
\end{equation}
Thus, the Fourier transform $U$ intertwines the actions of $\pi(M)'$ on $\mathcal{F}(E)\otimes_\pi H$ and on $\mathcal{F}(E^{\pi})\otimes_\iota H$ respectively.

The following theorem identifies the commutant of the Hardy algebra  represented by an induced representation.

\begin{thm}(~\cite{MuS3}, Theorem 3.9)\label{commutant of ind. rep.}
Let $E$ be a $W^*$-correspondence over $M$, and let $\pi:M\rightarrow B(H)$ be a faithful
normal representation of $M$ on a Hilbert space $H$.
Write $\rho_\pi$ for the representation $\pi^{\mathcal{F}(E)}$ of $H^{\infty}(E)$ on $\mathcal{F}(E)\otimes_{\pi}H$ induced by $\pi$, and write
$\rho^\pi$ for the representation of $H^{\infty}(E^{\pi})$ on
$\mathcal{F}(E)\otimes_\pi H$ defined by
\begin{equation}\label{Formula for rho^sigma}
\rho^\pi(X)=U^*\iota^{\mathcal{F}(E^{\pi})}(X)U,
\end{equation}
with $X\in H^{\infty}(E^{\pi})$.
Then $\rho^\pi$ is an
ultraweakly continuous, completely isometric representation of $H^{\infty}(E^{\pi})$ that extends the representation $\nu\times\Psi$ of $\mathcal{T}_+(E^{\pi})$, and
$\rho^\pi(H^{\infty}(E^{\pi}))$ is the commutant of $\rho_\pi(H^{\infty}(E))$, i.e. $\rho^\pi(H^{\infty}(E^{\pi}))=\rho_\pi(H^{\infty}(E))'$.
\end{thm}

\begin{cor}(~\cite{MuS3}, Corollary 3.10)\label{bicommutant of ind. rep.}
In the preceding notation,
$\rho_\pi(H^{\infty}(E))''=\rho_\pi(H^{\infty}(E))$.
\end{cor}

Combining this corollary with the well known fact that the commutant $\mathcal{A}'$ of every operator algebra $\mathcal{A}$ is WOT-closed, we obtain that $\rho_\pi(H^{\infty}(E))$ is WOT-closed.


\section{Reflexivity of the Hardy algebras.}
In this section we consider the reflexivity of the Hardy algebra $\rho_\pi(H^{\infty}(E))$. First we introduce the notions of a quantitative analog of reflexivity, the hyperreflexivity for operator algebras, and obtain some elementary consequences in our setting. After this observations we present our main results.

\subsection{Reflexivity for operator algebras.}

Let $H_0$ be a Hilbert and let $\mathcal{A}_0$ be some operator algebra, acting on it.
We write $H_0^{(\infty)}$ for the direct sum  of countable number of copies of a Hilbert space $H_0$. It can be naturally identified with the tensor product $H_0\otimes l^2$, were $l^2=l^2(\mathbb{Z}_+)$. The operator $A\otimes I\in \mathcal{A}_0\otimes I$ can be viewed as the infinite ampliation of $A$, ($A\in \mathcal{A}_0$). Then the algebra $\mathcal{A}_0\otimes I_{l^2}$ is unitarily equivalent to the infinite ampliation $\mathcal{A}_0^{(\infty)}$.
In ~\cite{RadRos1} Radjavi and Rosenthal proved that if the algebra $\mathcal{A}_0$ is unital and WOT-closed then the algebra $\mathcal{A}_0\otimes I$ is reflexive. In fact, this theorem may be strengthened by requiring that $\mathcal{A}_0$ will be only ultraweakly closed.

The operator algebra of the form $\mathcal{A}\otimes I$, where $I$ is the identity operator on an infinite dimensional Hilbert space $K$, is called an algebra of infinite multiplicity. Thus every WOT-closed or even ultraweakly closed operator algebra of infinite multiplicity is reflexive.

Returning to our general setting, let $M$ be a general $W^*$-algebra, $E$ be an arbitrary $W^*$ -correspondence over $M$ and $\pi$ be a faithful normal representation of $M$ on a Hilbert space $H$. Recall that we showed in previous section that the algebra $\rho_\pi(H^{\infty}(E))$ is WOT-closed. Assume that $\pi$ has an infinite multiplicity. This means that there is a separable Hilbert space $K$ such that $H=H_0\otimes K$ and there is a normal representation $\pi_0:M\rightarrow H_0$ such that $\pi=\pi_0\otimes I_K$. Let us show that in this case the algebra $\rho_\pi(H^{(\infty)}(E))$ has an infinite multiplicity. The space $\mathcal{F}(E)\otimes_\pi H$ may be written as $(\mathcal{F}(E)\otimes_{\pi_0}H_0)\otimes K$. Then the induced representation $\pi^{\mathcal{F}(E)}$ of the von Neumann algebra $\mathcal{L}(\mathcal{F}(E))$
has the form $(\pi_0\otimes I_K)^{\mathcal{F}(E)}:Z\mapsto Z\otimes I_H=
(Z\otimes I_{H_0})\otimes I_K=\pi_0^{\mathcal{F}(E)}(Z)\otimes I_K$. Hence $\mathcal{F}(E)\otimes_{\pi}H$ and $(\mathcal{F}(E)\otimes_{\pi_0}H_0)\otimes K$ are identified as $\mathcal{L}(\mathcal{F}(E))$-modules and $\pi^{\mathcal{F}(E)}(\cdot)=\pi_0^{\mathcal{F}(E)}(\cdot)\otimes I_K$ has an infinite multiplicity.
Conversely, for every faithful normal representation $\pi_0$ of $M$ on a Hilbert space $H_0$, the infinite ampliation $\rho_{\pi_0}(H^{\infty}(E))^{(\infty)}$ is acting on the Hilbert space $H_0^{(\infty)}\cong H_0\otimes l^2$, and may be identified with the algebra $\rho_{\pi_0}(H^{\infty}(E))\otimes I_{l^2}$.
Since the algebra $\rho_{\pi_0}(H^{\infty}(E))$ is WOT-closed, we get that $\rho_\pi(H^{\infty}(E))=H^{\infty}(E)\otimes I_H=(H^{\infty}(E)\otimes I_{H_0})\otimes I_K$ is reflexive.

In fact, every unital WOT-closed operator algebra of infinite multiplicity is hyperreflexive. This notion is a quantitative version of reflexivity.

Let $\mathcal{A}\subseteq B(H)$ be an operator algebra acting on a Hilbert space $H$.
For every $T\in B(H)$ set
$$\beta(T,\mathcal{A}):=\sup\{\|P^{\perp}TP\|: P\in Lat\mathcal{A}\}.$$
Then, for every $T\in B(H)$, $\beta(T,\mathcal{A})\leq dist(T,\mathcal{A})$, since if $P\in Lat\mathcal{A}$ then, for every $S\in \mathcal{A}$, $\|P^\perp TP\|=\|P^\perp (T-S)P\|\leq \|T-S\|$. Then, for every $S\in \mathcal{A}$, $\beta(T,\mathcal{A})\leq \|T-S\|$. Clearly,
$\beta(T,\mathcal{A})$ defines a seminorm on $B(H)$. Clearly also, $\beta(T,\mathcal{A})=0$ if and only if $T\in Alg\ Lat\ \mathcal{A}$.

In these terms we may redefine the reflexivity of the algebra $\mathcal{A}$ as follows:

\begin{equation*}
\mathcal{A}\ \text{is reflexive if and only if}\ \mathcal{A}=\{ T\in B(H): \beta(T,\mathcal{A})=0 \}
\end{equation*}
By definition, the algebra $\mathcal{A}$ is hyperreflexive if there is a constant $C\geq 0$, which is independent of $T$, such that
\begin{equation*}
\beta(T,\mathcal{A})\leq dist(\mathcal{A}, T)\leq C\beta(T,\mathcal{A}).
\end{equation*}
The infimum $C_T$ of all such numbers $C$ is called the constant of hyperreflexivity or the distance constant.
Clearly, each hyperreflexive algebra is reflexive. It is known that the following classes of algebras are hyperreflexive: the nest algebras (Arveson, ~\cite{Arv2}), the algebra
$\mathcal{L}_n$ (Davidson in ~\cite{Dav2} for $n=1$ and Davidson and Pitts in ~\cite{DavP2} for $n>1$). M. Kennedy in ~\cite{Ken} proved hyperreflexivity of  some class of free semigroup algebras. In ~\cite{KriPow} D. Kribs and S. Power show the hyperreflexivity of $\mathcal{L}_G$ for some special case of graphs, and in ~\cite{JaePow} F. Jaeck and S. Power shows the hyperreflexivity of $\mathcal{L}_G$ for any finite graph $G$.
Note also that the problem of characterization of von Neumann algebras which are hyperreflexive still open and equivalent to numerous long standing unsolved classical problems. But some partial results are known, for example, every injective von Neumann algebra is hyperreflexive.

The fact that every WOT-closed algebra $\mathcal{A}$ of infinite multiplicity is hyperreflexive is folklore (see ~\cite{DavP2}), and can be found, for example, in ~\cite{DavP2} (with distance constant at most 9), and with a a very short proof in ~\cite{JaePow} (with distance constant 3).

\begin{thm}\label{Hyperrefl of HardyAlg for infinitemultipl of pi} Let $\pi:M\rightarrow B(H)$ be a faithful normal representation of infinite multiplicity. Then the algebra
$\rho_\pi(H^{\infty}(E))$ is hyperreflexive with distance constant at most 3.
\end{thm}
\noindent{Proof.} Indeed, we saw $\rho_\pi(H^\infty(E))$ is WOT-closed and in the assumptions of the theorem, our algebra $\rho_\pi(H^\infty(E))$ has an infinite multiplicity. Hence it is hypereflexive with the distance constant 3.
\begin{flushright}
$\square$
\end{flushright}

Let us consider the special case when the algebra $M$ is a factor of type $III$. Since $M$ is a factor, then for every two projections  $p$ and $q$ in $M$ one has either $p\preceq q$ or $q\preceq p$ (in the sense of Murray-von Neumann) and since $M$ is of type $III$ all nonzero projections in $M$ are infinite and equivalent. We shall use the fact that every nonzero projection $p$ in such a factor can be ``divided by $\aleph_0$", i.e. there is a sequence $\{p_i\}_{i\geq 1}\subset M$ of pairwise orthogonal subprojections of $p$, such that $p=\sum p_i$ and $p_i\sim p$ for every $i$.

\begin{cor}
Let $M$ be a factor of type III. Then the algebra $\rho_\pi(H^{\infty}(E))$ is hyperreflexive.
\end{cor}
\noindent{Proof.}
Since $M$ is a factor of type $III$, so is $\pi(M)'$. It follows that $\pi$ is of infinite multiplicity.
To see this, write $I=\sum_{i\geq 0}p_i$ where $\{p_i\}_{i=0}^{\infty}$ are pairwise orthogonal, equivalent projections in $\pi(M)'$. Thus, there are partial isometries $\{u_j\}_{j=0}^{\infty}$ in $\pi(M)'$ with $u_j^*u_j=p_0$ and $u_ju_j^*=p_j$. Writing $H_0$ for $p_0H$ and $\pi_0$ for $\pi|_{H_0}$ we easily see that $\pi$ is unitarily equivalent to $\pi_0\otimes I_{l^2}$. Hence, $\pi$ is of infinite multiplicity. It follows from Theorem $\ref{Hyperrefl of HardyAlg for infinitemultipl of pi}$ that $\rho_\pi(H^{\infty}(E))$ is hyperreflexive and the distance constant is at most 3.

\begin{flushright}
$\square$
\end{flushright}

\subsection{Main results.}
Consider the algebra $Alg\ Lat\ \rho(H^{\infty}(E))$ where $\rho=\rho_\pi$ is an induced representation of $H^{\infty}(E)$ defined by the faithful normal representation $\pi$.

First we shall show that every $Z\in Alg\ Lat\ \rho(H^{\infty}(E))$ lies in $\rho(\mathcal{L}(\mathcal{F}(E)))$, hence, $Z$ has the form $Z=Y\otimes I_H$, where $Y$ is some element in $\mathcal{L}(\mathcal{F}(E))$.
We need the following simple auxiliary lemma.
\begin{lem}
Let $B$ be a von Neumann algebra acting on a Hilbert space $H$, and
let $\mathcal{A}\subset B(H)$ be some operator algebra. Assume that
$B\subset \mathcal{A}'$. Then $Alg\ Lat\ \mathcal{A}\subset B'$.
\end{lem}

\noindent{Proof.} For each $b\in B$ and each $a\in \mathcal{A}$, $ba=ab$. In
particular, for each projection $p\in B$  we have $pa=ap$, and the range of such
projection is in $Lat\ \mathcal{A}$. Let $\mathcal{M}\in Lat\ \mathcal{A}$ and
let $p_{\mathcal{M}}$ be
its projection. Then for each $c\in Alg\ Lat\ \mathcal{A}$ we have
$cp_{\mathcal{M}}=p_{\mathcal{M}}cp_{\mathcal{M}}$. In particular
$(1-p)cp=0$ for every projection $p\in B$. Since $B$ is a von Neumann algebra we
have also $pc(1-p)=0$ for all $p\in B$. Hence $cp=pcp=pc$. So, every $c\in Alg\ Lat\ \mathcal{A}$
commutes with every projection $p\in B$. It follows that $c\in B'$.
\begin{flushright}
$\square$
\end{flushright}

Let $\iota$ be a identity representation of $\pi(M)'$ on the Hilbert space $H$. Recall from the Preliminary section that for $b\in \pi(M)'$ the formula
$b\mapsto I_{\mathcal{F}(E)}\otimes \iota(b)$ defines a faithful normal representation
of $\pi(M)'$ on $\mathcal{F}(E)\otimes_\pi H$. We write $I_{\infty}\otimes \iota$ for this representation. Similarly, we write $I_n\otimes \iota$ for subsrepresentation of $\pi(M)'$ on
$E^{\otimes n}\otimes_{\pi}H$. Thus, if $\xi\in E^{\otimes n}$, then $(I_n\otimes \iota(b))(\xi\otimes h)=\xi\otimes bh$, $n\geq 0$. Frequently we shall drop the letter $\iota$, writing $I_n\otimes b$ for $I_n\otimes \iota(b)$.

\begin{cor}\label{Alg Lat in commutant of M'}
$Alg\ Lat\ \rho(H^{\infty}(E))\subset (I_{\infty}\otimes
\pi(M)')'=\mathcal{L}(\mathcal{F}(E))\otimes I_H$.

In particular, every $Z\in Alg\ Lat\ \rho_\pi(H^{\infty}(E))$ is of the form
$Y\otimes I_H$, for some $Y\in \mathcal{L}(\mathcal{F}(E))$.
\end{cor}

\noindent{Proof.} Let us take in the previous lemma $B:=(I_{\infty}\otimes \pi(M)')$
and $\mathcal{A}:= \rho(H^{\infty}(E))=H^{\infty}(E)\otimes I_H$.
Clearly, $I_{\infty}\otimes \pi(M)'\subset \rho(H^{\infty}(E))'$.
Hence
$$Alg\ Lat\ \rho(H^{\infty}(E))\subset (I_{\infty}\otimes \pi(M)')'.$$
By Rieffel's Theorem $\ref{Rieffel thm}$ we have $(I_{\infty}\otimes
\pi(M)')'=\mathcal{L}(\mathcal{F}(E))\otimes I_H$. Thus,
$Alg\ Lat\ \rho_\pi(H^{\infty}(E))\subset \mathcal{L}(\mathcal{F}(E))\otimes I_H$,
and the corollary follows.
\begin{flushright}
$\square$
\end{flushright}
In our attempts to prove the reflexivity of $\rho(H^{\infty}(E))$ we try to generalize the proof of Arias and Popescu from ~\cite{APo}. Thus, with each
$Y\otimes I_H\in Alg\ Lat\ \rho(H^{\infty}(E))$ we associate the series of its Fourier coefficients.

Namely, according to ~\cite{MuS3}, with every $T\in
\mathcal{L}(\mathcal{F}(E))$ we associate the series of operators $\Phi_j(T)$, that will be called Fourier coefficients, as follows.

Let $\{W_t: t\in \mathbb{R}\}$ be the one-parameter unitary group in
$\mathcal{L}(\mathcal{F}(E))$, defined by
$$W_t:=\sum_{n=0}^{\infty}e^{int}P_n.$$
Here $P_n$ is the projection on the n-th summand in $\mathcal{F}(E)$. One can
check that this series converges in the $w^*$-topology in
$\mathcal{L}(\mathcal{F}(E))$.
Further, let $\gamma_t(Y)=Ad W_t(Y)=W_tYW_t^*$. Then $\{\gamma_t:t\in \mathbb{R}\}$, is
a $w^*$-continuous action of $\mathbb{R}$ on $\mathcal{L}(\mathcal{F}(E))$,
called the gauge automorphism group.

The $j$-th Fourier coefficient $\Phi_j(T)$ is defined by
$$\Phi_j(T)=(1/2\pi)\int_0^{2\pi}e^{-ijt}\gamma_t(T)dt.$$
where the integral converges in the $w^*$-topology in
$\mathcal{L}(\mathcal{F}(E))$. Simple calculation gives us

\begin{equation}\label{Phi_j(T)=sum(P_k+jTP_k$)}
\Phi_j(T)=\sum_kP_{k+j}TP_k.
\end{equation}

\begin{lem}\label{Phi_j in Alg Lat}
Let $Y\otimes I_H\in Alg\ Lat\ \rho(H^{\infty}(E))$, then for each $j=0,1,2,...$,
the operator $\Phi_j(Y)\otimes I_H$ is in $Alg\ Lat\ \rho(H^{\infty}(E))$.
\end{lem}
\noindent{Proof.} A direct calculation gives the identity
$W_tT_{\xi}=e^{it}T_{\xi}W_t$, for $\xi\in E$, $t\in \mathbb{R}$. Note that for
every $t$, $W_t^*=W_{-t}$ and that $W_t$ has a closed range . It follows that if
$\mathcal{M}\in Lat\ \rho(H^{\infty}(E))$ then also
$\mathcal{M}_t:=(W_t\otimes I_H)\mathcal{M}$ is in $Lat\ \rho(H^{\infty}(E))$.

Let $Y\otimes I_H\in Alg\ Lat\ \rho(H^{\infty}(E))$.
Clearly, it is enough to show that also $\gamma_t(Y)\otimes I_H\in Alg\ Lat\
\rho(H^{\infty}(E))$.
We have $(W_tYW_t^*\otimes I_H)\mathcal{M}=(W_tY\otimes I_H)\mathcal{M}_{-t}
\subseteq (W_t\otimes I_H)\mathcal{M}_{-t}=(W_t\otimes I_H)(W_t^*\otimes
I_H)\mathcal{M}=\mathcal{M}$. So,
$(\gamma_t(Y)\otimes I_H)\mathcal{M}\subseteq \mathcal{M}$.
\begin{flushright}
$\square$
\end{flushright}
The operator $\Phi_j(Y)\otimes I_H$ will be called the $j$-th Fourier coefficient of $Y\otimes I_H\in Alg\ Lat\ \rho(H^{\infty}(E))$. Note that if $Y\otimes I_H\in Alg\ Lat\ \rho(H^{\infty}(E))$, then from the
formula $\Phi_j(Y)=\sum_kP_{k+j}YP_k$ it follows that $\Phi_j(Y)=0$ for every
$j<0$.

To prove the reflexivity of $\rho(H^{\infty}(E))$ it is enough to show that, given $Y\otimes
I_H\in Alg\ Lat\ \rho(H^{\infty}(E))$, then every $\Phi_j(Y)\otimes I_H$,
$j=0,1,2,...$, is in $\rho(H^{\infty}(E))$.
\begin{lem}\label{Appl of Fejer}
Let $Y\otimes I_H\in Alg\ Lat\ H^{\infty}(E)\otimes I_H$ and assume that
each $\Phi_j(Y)\otimes I_H\in \rho(H^{\infty}(E))$, $j\geq 0$. Then
$Y\otimes I_H\in \rho(H^{\infty}(E))$.
\end{lem}
\noindent{Proof.} It follows from
~\cite[p.366]{MuS3}, that the $k$-th arithmetic mean operators $\sigma_k(Y):=\sum_{|j|<k}(1-\frac{|j|}{k})\Phi_j(Y)$, $k\geq 0$, tend to $Y$ in the weak$^*$-topology in $\mathcal{L}(\mathcal{F}(E))$ as $k\rightarrow \infty$. But then the operators $\sigma_k(Y)\otimes I_H$ tend ultraweakly to $Y\otimes I_H\in B(\mathcal{F}(E)\otimes H)$. Hence, if all $\Phi_j(Y)\otimes I_H$ are in $\rho(H^{\infty}(E))$ then so is $Y\otimes I_H$.
\begin{flushright}
$\square$
\end{flushright}
Let $\xi\in E$, and let $L_{\xi}$ be the operator on $H$ defined by
$L_{\xi}:h\mapsto \xi\otimes h$. By $L^{(k)}_{\xi}$ we denote the operator
$L^{(k)}_{\xi}:E^{\otimes k}\otimes_{\pi}H\rightarrow E^{\otimes
k+1}\otimes_{\pi}H$, that on the elementary tensors is defined by $\eta\otimes
h\mapsto \xi\otimes\eta\otimes h$. So, $L_{\xi}=L^{(0)}_{\xi}$. Analogously, for
arbitrary $x\in E^{\otimes k}$, we let $L_x^{(n)}$ denote the operator from
$E^{\otimes n}\otimes_\pi H$ to $E^{\otimes n+k}\otimes_\pi H$ defined by
$\zeta\otimes h\mapsto x\otimes(\zeta\otimes h)$, where $\zeta\in E^{\otimes n}$.
For simplicity we often write $L^{(0)}_x=L_x$.

In the following theorem, which is a main result in this subsection, we agree to write $\xi^{(n)}$ for a general element of the correspondence $E^{\otimes n}$, when $n\geq 2$.

\begin{thm}\label{Repres. of Phi_j} Let $Y\otimes I_H\in Alg\ Lat\
\rho(H^{\infty}(E))$. Then

1) For $j=0$ there exists a sequence $\{a_s\}_{s\geq 0}\subset M$ such that
$(\Phi_0(Y)\otimes I_H)|_H=\pi(a_0)$ for $s=0$, and $(\Phi_0(Y)\otimes I_H)|_{E^{\otimes s}\otimes_{\pi}H}=\phi_s(a_s)\otimes I_H$ for every $s>0$,.
Hence, with respect to the decomposition
$\mathcal{F}(E)\otimes_{\pi}H=\sum^{\oplus}_{s\geq 0}E^{\otimes s}\otimes_{\pi}H$
the operator
$\Phi_0(Y)\otimes I_H$ is represented by the diagonal matrix
\begin{equation}\label{Matrix repres for Phi_0}
\text{diag}(\pi(a_0),\phi_1(a_1)\otimes I_H,...).
\end{equation}

2) For $j\geq 1$ there exists a sequence $\{\xi^{(j)}_s\}_{s\geq 0}\subset
E^{\otimes j}$ such that
$(\Phi_j(Y)\otimes I_H)|_{E^{\otimes
s}\otimes_{\pi}H}=(T_{\xi^{(j)}_s}\otimes I_H)|_{E^{\otimes
s}\otimes_{\pi}H}=L^{(s)}_{\xi^{(j)}_s}$.
Thus, with respect to the decomposition
$\mathcal{F}(E)\otimes_{\pi}H=\sum^{\oplus}_{s\geq 0}E^{\otimes
s}\otimes_{\pi}H$, the operator
$\Phi_j(Y)\otimes I_H$ is represented by the "j"-th subdiagonal matrix

\begin{equation}\label{Matrix repres for Phi_j}
\left(
  \begin{array}{cccccc}
    0   & 0   & 0   & 0   & ... & ... \\
    ... & ... & ... & ... & ... & ... \\
    0   & ... & ... & ... & ... & ... \\
    L^{(0)}_{\xi^{(j)}_0}   & 0   & 0   & 0 & ... & ... \\
    0   & L^{(1)}_{\xi^{(j)}_1}   & 0   & 0 & ... & ... \\
    0   & 0   & L^{(2)}_{\xi^{(j)}_2}   & 0 & ... & ... \\
    ... & ... & ... & ... & ... & ... \\
  \end{array}
\right)
\end{equation}
\end{thm}

\noindent{Proof.} Recall first from Lemma $\ref{Phi_j in Alg Lat}$ that if $Y\otimes I_H\in Alg\ Lat\ \rho(H^\infty)$ so is $\Phi_j(Y)\otimes I_H$ for every $j\geq 0$.
(1)  Since $\Phi_0(Y)=\sum_kP_kYP_k$, we obtain that $(\Phi_0(Y)\otimes
I_H)|_{E^{\otimes k}\otimes_{\pi}H}=P_kYP_k$, i.e. each summand $E^{\otimes
k}\otimes_{\pi}H$ of $\mathcal{F}(E)\otimes_{\pi}H$ is $\Phi_0(Y)\otimes I_H$-
invariant.
Consider the restriction $(\Phi_0(Y)\otimes I_H)|_H$ where $H\cong M\otimes_\pi H$. The representation
$I_{\infty}\otimes \iota$ of $\pi(M)'$ restricted to $H$ is simply $\iota$, the
identity representation of $\pi(M)'$ on $H$. Hence, by Corollary $\ref{Alg Lat in commutant of M'}$ and since $\Phi_j(Y)\otimes I_H\in Alg\ Lat\ \rho(H^\infty)$, we get the identity
$(\Phi_0(Y)\otimes
I_H)(I_{\infty}\otimes\iota(b))=(I_{\infty}\otimes\iota(b))(\Phi_0(Y)\otimes
I_H),\, b\in \pi(M)'$, which when
restricted to $H$, yields
$$(\Phi_0(Y)\otimes I_H)b=b(\Phi_0(Y)\otimes I_H),\,\, b\in \pi(M)'.$$
Thus, $(\Phi_0(Y)\otimes I_H)|_H=\pi(a_0)$, for some $a_0\in M$.

Let $n\geq 1$. For the restriction $(\Phi_0(Y)\otimes I_H)|_{E^{\otimes n}\otimes_{\pi}H}$ we have
$$(\Phi_0(Y)\otimes I_H)(I_{\infty}\otimes b)|_{E^{\otimes
n}\otimes_{\pi}H}=(P_nYP_n\otimes I_H)(I_n\otimes b)=(I_n\otimes b)(P_nYP_n\otimes
I_H),$$
where $I_n$ is the identity operator $I_{E^{\otimes n}}$.
We write $S_n\otimes I_H$ for the restriction $(\Phi_0(Y)\otimes I_H)|_{E^{\otimes
n}\otimes_{\pi}H}$, where $S_n=P_nYP_n\in \mathcal{L}(E^{\otimes n})$.

We want to show that for every $n\geq 1$, $S_n=\phi_n(a_n)$ for some $a_n\in M$.
Take $x\in E^{\otimes n}\otimes_\pi H$ be arbitrary and set
$$\mathcal{M}_{x}:=\overline{\rho(H^{\infty}(E))x}.$$

So, $\mathcal{M}_{x}$ is a $\rho(H^{\infty}(E))$-invariant subspace in $\mathcal{F}(E)\otimes_{\pi}H$ and can be
written as
$$\mathcal{M}_{x}=\overline{(\phi_n(M)\otimes I_H)x}\oplus \overline{E\otimes_{\phi_n\otimes I_H}x}\oplus...$$

We see that
$\mathcal{M}_x\cap(E^{\otimes n}\otimes_{\pi}H)=\overline{(\phi_n(M)\otimes I_H)x}$
and is invariant under $\Phi_0(Y)\otimes I_H$ because $\Phi_0\otimes I_H\in Alg\ Lat\ \rho(H^\infty(E))$.
Thus,
$(S_n\otimes I_H)(x)\in \overline{(\phi_n(M)\otimes
I_H)x}$ and we obtain that there is a net
$(a_\alpha)\subset M$ such that

$$(S_n\otimes I_H)x=\lim_{\alpha}(\phi_n(a_\alpha)\otimes I_H)x,$$
and the net $(a_{\alpha})$ depends on the choice of $x\in
E^{\otimes n}\otimes_\pi H$.

Fix any projection $p\in (\phi_n(M)\otimes I_H)'$. Replacing $x$ by $px$ and by $p^{\perp}x$ (where $p^{\perp}=I-p$), we get two nets $(b_{\alpha})$ and $(c_{\alpha})$ in $M$ such that
$$(S_n\otimes I_H)(px)=\lim_\alpha (\phi_n(b_\alpha)\otimes
I_H)(px),$$
$$(S_n\otimes I_H)(p^\perp x)=\lim_\alpha (\phi_n(c_\alpha)\otimes
I_H)(p^\perp x),$$

Then
$$(S_n\otimes I_H)x=\lim_\alpha (\phi_n(b_\alpha)\otimes I_H)(px)+\lim_\alpha
\phi_n(c_\alpha)\otimes I_H)(p^\perp x).$$

Now applying $p$ on both sides and using the facts that $p\in (\phi_n(M)\otimes
I_H)'$ and $pp^{\perp}=0$ we obtain
$$p(S_n\otimes I_H)x=\lim_\alpha (\phi_n(b_\alpha)\otimes
I_H)(px)=(S_n\otimes I_H)px.$$
Hence for every projection $p\in (\phi_n(M)\otimes I_H)'$ we have $p(S_n\otimes
I_H)=(S_n\otimes I_H)p$ at $x\in E^{\otimes n}\otimes_{\pi}H$. Since the choice of $x$ is arbitrary we get
$$p(S_n\otimes I_H)=(S_n\otimes I_H)p,$$
on $E^{\otimes n}\otimes_{\pi}H$. Since $p\in (\phi_n(M)\otimes I_H)'$ is an arbitrary projection we get:
$$S_n\otimes I_H\in (\phi_n(M)\otimes I_H)''=\phi_n(M)\otimes I_H,$$
and this implies that there exists $a_n\in M$ such that
$$S_n\otimes I_H=\phi(a_n)\otimes I_H.$$
We proved that $(\Phi_0\otimes I_H)|_{E^{\otimes n}\otimes_\pi H}=\phi_n(a_n)\otimes I_H$ for some $a_n\in M$.

Hence, $\Phi_0(Y)=diag(\pi(a_0), \phi_1(a_1)\otimes I_H,\phi_2(a_2)\otimes I_H,...)$.
\vskip 0.5 cm
(2) Let $j\geq 1$. From the formula $\Phi_j(Y)=\sum_{k}P_{k+j}YP_{k}$ we see that
for every $s\geq0$
$$\Phi_j(Y)\otimes I_H:E^{\otimes s}\otimes_\pi H\rightarrow E^{\otimes s+j}\otimes_\pi H.$$
Consider the restriction $(\Phi_j(Y)\otimes I_H)|_H$.
Thus,
$\Phi_j(Y)\otimes I_H|_H$ acts from $H$ to
$E^{\otimes j}\otimes_{\pi}H$.

By Corollary $\ref{Alg Lat in commutant of M'}$ and Lemma $\ref{Phi_j in Alg
Lat}$ we have
$$(\Phi_j(Y)\otimes I_H)(I_{\infty}\otimes b)x=(I_{\infty}\otimes
b)(\Phi_j(Y)\otimes I_H)x,$$
where $x\in \mathcal{F}(E)\otimes_{\pi}H$ and $b\in \pi(M)'$.

In particular, for every $h\in H$

$$(\Phi_j(Y)\otimes I_H)(bh)=(I_j\otimes b)(\Phi_j(Y)\otimes I_H)h.$$

Let $U=U_\pi:\mathcal{F}(E)\otimes_{\pi}H\rightarrow
\mathcal{F}(E^{\pi})\otimes_{\iota}H$
be the Fourier transform defined by $\pi$ (see Definition $\ref{Fourier trans definition}$). It is pointed out in Chapter 2 that
$U$ is a Hilbert space isomorphism as well as each
restriction
$U_s:=U|_{E^{\otimes s}\otimes_{\pi}H}:E^{\otimes s}\otimes_{\pi}H\rightarrow
(E^{\pi})^{\otimes s}\otimes_{\iota}H$. From formula $(\ref{Intertwining M' and U})$, the operator $U$ intertwines the representations $I_{\infty}\otimes
\iota$
and $\iota^{\mathcal{F}(E^{\pi})}$ of $\pi(M)'$ on $\mathcal{F}(E)\otimes_{\pi}H$
and
$\mathcal{F}(E^{\pi})\otimes_{\iota}H$ respectively.

Regarding $H$ and $(E^{\pi})^{\otimes j}\otimes_{\iota}H$ as
subspaces in $\mathcal{F}(E)\otimes_{\pi}H$ and
$\mathcal{F}(E^{\pi})\otimes_{\iota}H$ respectively,
we consider  the operator
$$U(\Phi_j(Y)\otimes I_H)|_H:H\rightarrow (E^{\pi})^{\otimes
j}\otimes_{\iota}H.$$

Hence, $U(\Phi_j(Y)\otimes I_H)|_H$ intertwines the actions $\iota$ and
$\iota^{(E^{\pi})^{\otimes j}}$ of $\pi(M)'$.
It follows that
$U(\Phi_j(Y)\otimes I_H)|_H$ is contained in the second dual $(E^{\otimes j})^{\pi,\iota}$ and, by the duality theory $(E^{\otimes j})^{\pi,\iota}\cong E^{\otimes j}$. Hence, there exists a unique
$\xi^{(j)}\in E^{\otimes j}$, that corresponds
to $U(\Phi_j(Y)\otimes I_H)|_H$. We write
$\widehat{\xi^{(j)}}=U(\Phi_j(Y)\otimes I_H)|_H$.

We will show that $(\Phi_j(Y)\otimes I_H)|_H=(T_{\xi^{(j)}}\otimes I_H)|_H=L_{\xi^{(j)}}$. To
this end recall that by formula $(\ref{Formula for hat(xi)})$ form Chapter 2, we have
$$\widehat{\xi^{(j)}}=U_j\circ
L_{\xi^{(j)}}.$$
Thus $\widehat{\xi^{(j)}}=U_j(\Phi_j(Y)\otimes I_H)|_H=U_j\circ
L_{\xi^{(j)}}$, and we get

\begin{equation}\label{Restriction of Phi_j to H}
L_{\xi^{(j)}}^{(0)}=(T_{\xi^{(j)}}\otimes I_H)|_H=(\Phi_j(Y)\otimes I_H)|_H.
\end{equation}

We write $\xi^{(j)}_0$ for the $\xi^{(j)}$ obtained above.

Now we consider the restriction $S_n\otimes I_H:=(\Phi_j(Y)\otimes
I_H)|_{E^{\otimes n}\otimes_{\pi}H}$ for $n\geq 0$. The operator $S_n\otimes I_H=P_{n+j}YP_n\otimes I_H$ acts from $E^{\otimes n}\otimes_{\pi}H$ into $E^{\otimes j+n}\otimes_{\pi}H$.
Put $K_n:=E^{\otimes n}\otimes_{\pi}H$ and recall that $E^{\otimes j}\otimes_{\phi_n\otimes I_H}K_n=E^{\otimes
j}\otimes_{\phi_n\otimes I_H}(E^{\otimes n}\otimes_{\pi}H)\cong E^{\otimes
j+n}\otimes_{\pi}H$. Thus, $S_n\otimes I_H$ acts from $K_n$ into
$E^{\otimes
j}\otimes_{\phi_n\otimes I_H}K_n$.

Take $x\in K_n$ arbitrary and form the cyclic
$\rho(H^{\infty}(E))$-invariant subspace
$\mathcal{M}_{x}=\overline{\rho(H^{\infty}(E))x}$. The subspace $\mathcal{M}_x$ has the representation
$$\mathcal{M}_x=\overline{(\phi_n(M)\otimes I_H)x}\oplus \overline{E\otimes_{\phi_n\otimes I_H}x}\oplus...\oplus \overline{E^{\otimes
j}\otimes_{\phi_n\otimes I_H}x}\oplus...$$
Since $\Phi_j\otimes I_H\in Alg\ Lat\ \rho(H^\infty(E))$ we get $(\Phi_j\otimes I_H)\mathcal{M}_x\subseteq \mathcal{M}_x$.
Then $(S_n\otimes I_H)x\in (E^{\otimes j}\otimes K_n)\cap \mathcal{M}_x=\overline{E^{\otimes
j}\otimes_{\phi_n\otimes I_H}x}$, and there exists some net
$(\theta_{\alpha}^{(j)})\subset E^{\otimes j}$ such that
$$(S_n\otimes I_H)x=\lim_{\alpha}\theta_{\alpha}^{(j)}\otimes
x=\lim_{\alpha}(T_{\theta_{\alpha}^{(j)}}\otimes I_H)x.$$

Let $p\in (\phi_n(M)\otimes I_H)'$ be any projection, then for every $\alpha$
$$(T_{\theta_{\alpha}^{(j)}}\otimes I_H)px=\theta_{\alpha}^{(j)}\otimes px=(I_j\otimes
p)(\theta_{\alpha}^{(j)}\otimes x)=(I_{j}\otimes
p)(T_{\theta_{\alpha}^{(j)}}\otimes I_H)x,$$
where $I_j$ denotes the identity operator
$I_{E^{\otimes j}}$.
Now take $px$ and $p^{\perp}x=(I-p)x$ in $K_n$, one can find a nets $(\zeta_{\alpha}^{(j)})$ and
$(\vartheta_{\alpha}^{(j)})$ in $E^{\otimes j}$ such that
$(S_n\otimes I_H)px=\lim_{\alpha}(T_{\zeta_{\alpha}^{(j)}}\otimes I_H)px$
and $(S_n\otimes I_H)p^{\perp}x=\lim_{\alpha}(T_{\vartheta_{\alpha}^{(j)}}\otimes I_H)p^{\perp}x$. So,
$$\lim_{\alpha}(T_{\theta_{\alpha}^{(j)}}\otimes I_H)x=\lim_{\alpha}(T_{\zeta_{\alpha}^{(j)}}\otimes I_H)px
+\lim_{\alpha}(T_{\vartheta_{\alpha}^{(j)}}\otimes I_H)p^{\perp}x.$$

Applying $I_{j}\otimes p$ on both sides, we obtain

$$(I_{j}\otimes p)(S_n\otimes I_H)x=(I_{j}\otimes p)\lim_{\alpha}(T_{\theta_{\alpha}^{(j)}}\otimes
I_H)x=$$
$$=\lim_{\alpha}(T_{\zeta_{\alpha}^{(j)}}\otimes I_H)px=(S_n\otimes
I_H)px.$$

Since $p$ is an arbitrary projection in $(\phi_n(M)\otimes I_H)'$ and this holds for
every $x\in K_n$, we obtain that $S_n\otimes I_H$ intertwines the action of
$(\phi_n(M)\otimes I_H)'$:

\begin{equation}\label{Intertwining S_ntensor I_H and (phi_n(M)tens I_h)'}
(I_j\otimes b)(S_n\otimes I_H)=(S_n\otimes I_H)b,\,\,\,\forall{b}\in
(\phi_n(M)\otimes I_H)'.
\end{equation}

Denote $\iota_n:=id(\phi_n(M)\otimes I_H)'$ the identity representation of $(\phi_n(M)\otimes I_H)'$ on $K_n$ and let
 $$U_{\phi_n\otimes I_H}:E^{\otimes j}\otimes_{\phi_n\otimes
I_H}K_n\rightarrow (E^{\phi_n\otimes I_H})^{\otimes j}\otimes_{\iota_n}K_n,$$
be the Fourier transform defined by $\phi_n(\cdot)\otimes I_H$. By formula $(\ref{Intertwining M' and U})$ from Preliminaries, $U_{\phi_n\otimes I_H}$
intertwines the actions of $(\phi_n(M)\otimes I_H)'$:
$$U_{\phi_n\otimes I_H}(I_j\otimes b)=(\phi_{E^{\phi_n\otimes I_H}}(b)\otimes I_H)U_{\phi_n\otimes I_H},$$
for every $b\in (\phi_n(M)\otimes I_H)'$.

Thus, the composition $U_{\phi_n\otimes I_H}(S_n\otimes I_H)$ is in the second dual
$(E^{\otimes j})^{\phi_n\otimes I_H,\iota_n}$. By duality there exists a unique
$\xi^{(j)}\in E^{\otimes j}$ such that

$$\widehat{\xi^{(j)}}^*(\tilde{\eta}_1\otimes...\otimes\tilde{\eta_j}\otimes k_n)=
\tilde{L}_{\xi^{(j)}}^*((I_{j-1}\otimes\tilde{\eta}_1)...(I_{2}\otimes\tilde{\eta}_{j-1})\tilde{\eta}_j(k_n)),$$
where $\tilde{L}_{\xi^{(j)}}=(T_{\xi^{(j)}}\otimes I_H)|_{K_n}$, $\tilde{\eta}\in
E^{\phi_n\otimes I_H}$ and $k_n\in K_n$. Thus,
$\widehat{\xi^{(j)}}=\tilde{U}_j\circ \tilde{L}_{\xi^{(j)}}$, where
we denote $\tilde{U}_j=U_{\phi_n\otimes I_H}|_{E^{\otimes j}\otimes_{\phi_n\otimes
I_H}K_n}$.
Hence $$\tilde{U}_j^*\widehat{\xi^{(j)}}=\tilde{L}_{\xi^{(j)}}.$$
We write $\xi^{(j)}_n$ for $\xi^{(j)}$ obtained above.

We obtain that
$$(\Phi_j(Y)\otimes I_H)|_{E^{\otimes n}\otimes_\pi H}=(T_{\xi^{(j)}_n}\otimes
I_H)|_{E^{\otimes n}\otimes_\pi H}=L^{(n)}_{\xi^{(j)}_n},$$
and this proves the matrix representation $(\ref{Matrix repres for Phi_j})$ and the proof of the theorem is complete.
\begin{flushright}
$\square$
\end{flushright}

Our next step is to consider the action of every $\Phi_j(Y)\otimes I_H$ and of its adjoint $(\Phi_j(Y)\otimes I_H)^*$ on a suitable
$\rho(H^{\infty}(E))$-coinvariant subspaces $\mathfrak{M}$. Observe that $(\Phi_j(Y)\otimes I_H)^*=\Phi_j(Y)^*\otimes I_H$.

Our considerations are based on the following simple facts. Let $\mathcal{A}$ be an operator algebra, acting in Hilbert space $H$. Then
$$Lat\ \mathcal{A}^*=(Lat\ \mathcal{A})^{\perp},$$
where as usual $\mathcal{A}^*$ is the algebra of adjoint of elements of $\mathcal{A}$ and $(Lat\ \mathcal{A})^{\perp}$ is the lattice of orthogonal complements of subspaces from $Lat\ \mathcal{A}$.

\begin{lem}\label{Perp of adjoint}
Let $\mathcal{A}$ be an operator algebra, acting in Hilbert space $H$.
Then $$(Alg\ Lat\ \mathcal{A})^*\subseteq Alg\ Lat\ \mathcal{A}^*.$$
\end{lem}
\noindent{Proof.}
Let $Y\in Alg\ Lat\ \mathcal{A}$. From $Lat\ \mathcal{A}^*=(Lat\ \mathcal{A})^{\perp}$ we get $(1-p_{\mathcal{M}^{\perp}})Y^*p_{\mathcal{M}^{\perp}}=0$, for every
$\mathcal{M}\in Lat\ \mathcal{A}$. So,
$Y^*\mathcal{M}^{\perp}\subseteq\mathcal{M}^{\perp}$. Finally, $Y^*\tilde{\mathcal{M}}\subseteq \tilde{\mathcal{M}}$ for every $Y\in Alg\ Lat
\mathcal{A}$ and every $\tilde{\mathcal{M}}\in Lat\ \mathcal{A}^*$.
\begin{flushright}
$\square$
\end{flushright}

For every $s\geq 0$, let
$H(s):=\bigvee\{\eta^{(s)}(h):\eta\in E^\pi, h\in H\},\ s=0,1,2,...$, where, clearly,
$H(0)=H$.
Put
\begin{equation}\label{DefnOfGenerSymmetricPart}
\mathfrak{M}=\sum_{s\geq 0}^{\oplus}H(s).
\end{equation}

Recalling the definition of the Fourier transform $U_{\pi}$ we see that
$H(s)=U_{\pi}^*\tilde{H}(s)$, where $\tilde{H}(s)=\bigvee\{\eta^{\otimes s}\otimes h:\eta\in E^{\pi},h\in H\}$. So, if we set $\widetilde{\mathfrak{M}}:=\sum^{\oplus}_{s\geq 0}\tilde{H}(s)$, then $\widetilde{\mathfrak{M}}=U_{\pi}(\mathfrak{M})$.

Notice that $H(0)=H$ and $H(1)=E\otimes_{\pi}H=U_{\pi}^*(\tilde{H}(s))=U_{\pi}^*(E^{\pi}\otimes_{\iota}H)$, which follows from ~\cite[Lemma 3.5]{MuS3}. But when $s\geq 2$ the subspace $H(s)$ in general is a proper subspace of $E^{\otimes s}\otimes_{\pi}H$.
Take for example $E=\mathbb{C}^n$ and $M=\mathbb{C}$. In this case $E^{\pi}$ is all bounded operators $B(H,\mathbb{C}^n\otimes_\pi H)$ and can be identified with the $n$-fold column space $C_n(B(H))$ over the algebra $B(H)$ (see ~\cite[Example 4.2]{MuS3}). If we take $H=\mathbb{C}$ then $C_n(B(\mathbb{C}))\cong \mathbb{C}^n$.
Hence, $(\mathbb{C}^{n})^{\otimes s}\otimes_\iota \mathbb{C}\cong (\mathbb{C}^{n})^{\otimes s}$ and the subspace $\tilde{H}(s)$ now is $\bigvee\{\eta\otimes...\otimes\eta:\eta\in \mathbb{C}^n\}$, that is, the symmetric tensor power of $\mathbb{C}^n$.  Thus, the subspace $\widetilde{\mathfrak{M}}$ is the ordinary symmetric Fock space of $\mathbb{C}^n$
and $\tilde{H}(s)\subsetneqq (\mathbb{C}^{n})^{\otimes s}$, for $s\geq 2$.

The example just described is the reason to call the subspace $\widetilde{\mathfrak{M}}$, as well as its Fourier image $\mathfrak{M}$, the symmetric part (subspace) of the corresponding full Fock space.

\begin{thm}\label{PropertOfSymmetricPart}
The subspace $\mathfrak{M}\subset \mathcal{F}(E)\otimes_\pi H$ is $\rho(H^{\infty}(E))$-coinvariant.
Let, further, $Y\otimes I_H\in Alg\ Lat\
\rho(H^{\infty}(E))$ and $\Phi_j(Y)\otimes I_H$ be it's $j$-th Fourier
coefficient. Then, in the notation of Theorem $\ref{Repres. of Phi_j}$,

a) if $j=0$ then
$$(\Phi_0(Y)^*\otimes I_H)|_{\mathfrak{M}}=
(\phi_{\infty}(a_0)^*\otimes I_H)|_{\mathfrak{M}};$$

b)  if $j\geq 1$, then
$$(\Phi_j(Y)^*\otimes I_H)|_{\mathfrak{M}}=
(T_{\xi^{(j)}_0}^*\otimes I_H)|_{\mathfrak{M}},$$
$\xi^{(j)}_0\in E^{\otimes j}$.
\end{thm}

For the proof of the theorem we shall need a couple of lemmas.
\begin{lem}\label{SimmetricMisCoinvariant} The subspace $\mathfrak{M}$ is $\rho(H^\infty(E))$-coinvariant.
\end{lem}
\noindent{Proof.} Let $a\in M$. For every $\eta\in E^{\pi}$ and every $s\geq 0$, one have
$$(\phi_\infty(a^*)\otimes I_H)\eta^{(s)}(h)=\eta^{(s)}(\pi(a^*)h).$$
Thus, every operator $\phi_\infty(a^*)\otimes I_H$ leaves invariant every subspace $H(s)$.

We shall use by the following simple formula (in fact, it was used in ~\cite{MuS3}).
Let $\eta\in E^{\pi}$ and $s\geq 1$ arbitrary. Then for every $j$ such that $1\leq j\leq s$ and every  $x\in E^{\otimes j}$,
\begin{equation}\label{General identity}
L_{x}^{(s-j)*}\eta^{(s)}(h)=\eta^{(s-j)}L_{x}^{(0)*}\eta^{(j)}(h),\,\,\,x\in
E^{\otimes j},\, 1\leq j\leq s,
\end{equation}
where $L_x^{(n)}$ are operators that were defined before Theorem $\ref{Repres. of Phi_j}$.
In particular, for $\xi\in E$,
$$L_{\xi}^*\eta^{(s)}(h)=\eta^{(s-1)}(L_{\xi}^{*}\eta(h)),\,\,\, s\geq 0.$$
Hence, for every $s\geq 0$, $\eta^{(s)}(h)\in H(s)$,
$$(T_\xi^*\otimes I_H)\eta^{(s)}(h)=\eta^{(s-1)}(L_{\xi}^{*}(h))\in H(s-1),$$
thus, for every $\xi\in E$, $(T_\xi^*\otimes I_H)H(s)\subset H(s-1)$, and this proves that $\mathfrak{M}$ is $\rho(H^{\infty}(E))$-coinvariant.
\begin{flushright}
$\square$
\end{flushright}

Now let us recall the definition of the Cauchy
transform from ~\cite{MuS7} (see also ~\cite{MuS3}).

Let $\eta\in E^{\pi}$, with $\|\eta\|<1$.
The Cauchy transform $C_\eta$ is the operator
from $H\cong M\otimes_{\pi}H$ into $\mathcal{F}(E)\otimes_{\pi}H$ defined by

$$C_{\eta}:h\mapsto h+\eta(h)+\eta^{(2)}(h)+\eta^{(3)}(h)+...,$$
where $\eta^{(k)}(h)=(I_{E^{\otimes k-1}}\otimes \eta)...(I_E\otimes
\eta)\eta(h)$.
It is pointed out in ~\cite{MuS7} that $C_{\eta}$ is bounded with $\|\tilde{C}_{\eta}\|\leq \frac{1}{1-\|\eta\|}$.
Thus, $C_{\eta}=column(I,\eta, \eta^{(2)}, \eta^{(3)},...)$.

It is easy to see that for every $k\geq 0$ and every $a\in M$ we have the equality
$$(I_{E^{\otimes k}}\otimes \eta)(\phi_k(a)\otimes I_H)=(\phi_{k+1}(a)\otimes
I_H)(I_{E^{\otimes k}}\otimes \eta).$$
Hence $(\phi_k(a)\otimes I_H)\eta^{(k)}(h)=\eta^{(k)}(\pi(a)h)$
for every $h\in H$, and we get
$$C_{\eta}\pi(a)=(\phi_{\infty}(a)\otimes I_H)C_{\eta},\,\, a\in M.$$
So, $C_\eta\in \mathcal{F}(E^\pi)$.
Let us show that $C_{\eta}$ has a closed range. To this end let
$K_{\eta}=(I_{\mathcal{F}(E)}\otimes \Delta_{*}(\eta))C_{\eta}$ be the Poisson kernel associated with $\eta$ as defined in ~\cite[Definition 8]{MuS7}. By Proposition 10 of that paper, $K_{\eta}$ is an isometry mapping from $H$ to
$\mathcal{F}(E)\otimes_{\pi}H$ and
$K_{\eta}^*K_{\eta}=C_{\eta}^*(I_{\mathcal{F}(E)}\otimes(\Delta_*(\eta))^2)C_{\eta}=I_H$.
So, $C_{\eta}^*(I_{\mathcal{F}(E)}\otimes(\Delta_*(\eta))^2)$ is the left inverse
of $C_{\eta}$, hence $C_\eta$ has a closed range.

For a given $\eta\in \mathbb{D}(E^{\pi})$ let
\begin{equation}\label{SubspaceM_eta}
\mathfrak{M}_{\eta}:=\bigvee\{C_{\eta}h:h\in H\}.
\end{equation}

Since $C_{\eta}$ has a closed range we get $\mathfrak{M}_{\eta}=C_{\eta}(H)$ and every element of $\mathfrak{M}_{\eta}$ has the form
$h+\eta(h)+\eta^{(2)}(h)+...$, for some $h\in H$.

Clearly, $\mathfrak{M}_\eta\subset \mathfrak{M}$ for every $\eta\in \mathbb{D}(E^\pi)$.
Moreover, each subspace $\mathfrak{M}_{\eta}$ is $\rho(H^{\infty}(E))$-coinvariant, which may be shown in the same way as  $\rho(H^{\infty}(E))$-coinvariance of $\mathfrak{M}$.

\begin{lem} \label{Action of Phi_j on M_eta} 1) Let $Y\otimes I_H\in Alg\ Lat\
\rho(H^{\infty}(E))$. Then, in the notations of Theorem $\ref{PropertOfSymmetricPart}$,

(i) for $j=0$,
$$(\Phi_0(Y)^*\otimes I_H)|_{\mathfrak{M}_{\eta}}=
(\phi_{\infty}(a_0)^*\otimes I_H)|_{\mathfrak{M}_{\eta}}.$$
and

(ii) for $j\geq 1$,
$$(\Phi_j(Y)^*\otimes I_H)|_{\mathfrak{M}_{\eta}}=
(T_{\xi^{(j)}_0}^*\otimes I_H)|_{\mathfrak{M}_{\eta}}.$$
\end{lem}
\noindent{Proof.} $\mathfrak{M}_{\eta}$ is $\rho(H^{\infty}(E))$-coinvariant,
and hence $(\Phi_j(Y)\otimes I_H)^*\mathfrak{M}_{\eta}\subseteq \mathfrak{M}_{\eta}$ for every $j\geq 0$. Take $x=\sum_{s\geq 0}\eta^{(s)}(h)\in
\mathfrak{M}_{\eta}$.

1) If $j=0$, then
$$(\Phi_0(Y)\otimes I_H)^*x=\sum_{s\geq 0}(\Phi_0(Y)^*\otimes I_H)\eta^{(s)}(h)=
\sum_{s\geq 0}\eta^{(s)}(k),$$
for some $k\in H$.
So, $$(\Phi_0(Y)^*\otimes I_H)\eta^{(s)}(h)=\eta^{(s)}(k),\,\, s\geq 0.$$
By part 1) of Theorem $\ref{Repres. of Phi_j}$, $(\Phi_{0}(Y)^*\otimes I_H)|H=\pi(a_0^*)$. Hence
$\pi(a_0^*)h=k$ and for every $s\geq 1$ we obtain, from the intertwining property of $\eta^{(s)}$,
$$\eta^{(s)}(k)=\eta^{(s)}(\pi(a_0^*)h)=(\phi_s(a_0^*)\otimes
I_H)\eta^{(s)}(h).$$
So, $$(\Phi_0(Y)^*\otimes I_H)\eta^{(s)}(h)=(\phi_s(a_0^*)\otimes
I_H)\eta^{(s)}(h),\,\, s\geq 1.$$

Hence,
$(\Phi_0(Y)\otimes
I_H)^*|_{\mathfrak{M}_{\eta}}=(\phi_{\infty}(a_0^*)\otimes
I_H)|_{\mathfrak{M}_{\eta}}$.

2) Let $j\geq 1$ and $x=\sum_{s\geq 0}\eta^{(s)}(h)\in\mathfrak{M}_{\eta}$ as in 1).
Since $\mathfrak{M}_{\eta}$, is $(\Phi_j(Y)\otimes I_H)^*$-invariant
we have
$$(\Phi_j(Y)^*\otimes I_H)x=\sum_{s\geq 0}\eta^{(s)}(k)\in\mathfrak{M}_{\eta},$$
for some $k\in H$.

From the inclusion $(\Phi_j(Y)\otimes I_H)(E^{\otimes s}\otimes_{\pi}H)\subseteq
E^{\otimes s+j}\otimes_{\pi}H$ we get $(\Phi_j(Y)^*\otimes I_H)\eta^{(l)}(h)=0$
for $l<j$.
Further,
$$(\Phi_j(Y)^*\otimes I_H)\eta^{(s+j)}(h)=\eta^{(s)}(k),\,\, s\geq 0,$$
and for $s=0$:
$$(\Phi_j(Y)^*\otimes I_H)\eta^{(j)}(h)=k.$$
By Theorem $\ref{Repres. of Phi_j}$, $2)$, $(\Phi_j(Y)^*\otimes
I_H)|_{E^{\otimes j}\otimes_\pi H}=(T^*_{\xi^{(j)}_0}\otimes I_H)|_{E^{\otimes
j}\otimes_\pi H}=L^{(0)*}_{\xi^{(j)}_0}$.
Hence,
$$k=L^{(0)*}_{\xi^{(j)}_0}\eta^{(j)}(h).$$
So, for $s\geq 0$
$$(\Phi_j(Y)^*\otimes I_H)\eta^{(s+j)}(h)=\eta^{(s)}(k)=
\eta^{(s)}(L^{(0)*}_{\xi^{(j)}_0}\eta^{(j)}(h)).$$
By identity $(\ref{General identity})$
$$(\Phi_j(Y)^*\otimes I_H)\eta^{(s+j)}(h)=
L^{(s)*}_{\xi^{(j)}_0}\eta^{(s+j)}(h).$$

Hence, for every $x\in \mathfrak{M}_{\eta}$,
$$(\Phi_j(Y)^*\otimes I_H)x=(T_{\xi^{(j)}_0}^*\otimes I_H)x,$$
and lemma follows.
\begin{flushright}
$\square$
\end{flushright}

\noindent{\textbf{Proof of Theorem $\ref{PropertOfSymmetricPart}$}.} By Lemma $\ref{SimmetricMisCoinvariant}$ $\mathfrak{M}$ is $\rho(H^\infty(E))$-coinvariant.
By the last lemma the restriction $(\Phi_j(Y)\otimes
I_H)^*|_{\mathfrak{M}_{\eta}}$ is of the form $(X^*\otimes I_H)|_{\mathfrak{M}_{\eta}}$, where $X$ is an element of $H^{\infty}(E)$ and is either of the form
$\phi_{\infty}(a)$ with $a\in M$,
or $T_{\xi^{(j)}}$ with $\xi^{(j)}\in E^{\otimes j}$, $j\geq 1$.
Note also that $X$ is independent of $\eta$ and, thus, this holds also for $\bigvee\{\mathfrak{M}_{\eta}:\eta\in E^\eta, \|\eta\|<1\}$. Clearly, this is also true for every subspace $H(s)$, $s=0,1,2,...$. Thus,
$(\Phi_j(Y)\otimes
I_H)^*|_{H(s)}=(X^*\otimes I_H)|_{H(s)},\,\, s\geq 0$.
We obtain
\begin{equation}\label{Action of Phi_j* on big M}
(\Phi_j(Y)\otimes
I_H)^*|_{\mathfrak{M}}=(X^*\otimes I_H)|_{\mathfrak{M}},
\end{equation}
where $X\in H^{\infty}(E)$ of the form pointed above.
\begin{flushright}
$\square$
\end{flushright}

\begin{prop}\label{Action of Phi_j on big M}
Let $Y\otimes I_H\in Alg\
Lat\ \rho(H^{\infty}(E))$ as above
and let $\Phi_j(Y)\otimes I_H$ be its $j$-th Fourier coefficient, represented as in
Theorem $\ref{Repres. of Phi_j}$. Let $Q=P_{\mathfrak{M}}$ the projection onto $\mathfrak{M}$.

1) If $j=0$ then $$Q(\Phi_0(Y)\otimes
I_H)|_{\mathfrak{M}}=(\phi_{\infty}(a_0)\otimes I_H)|_{\mathfrak{M}}.$$

2) If $j\geq 1$ then
$$Q(\Phi_j(Y)\otimes
I_H)|_{\mathfrak{M}}=Q(T_{\xi^{(j)}_0}\otimes I_H)|_{\mathfrak{M}}.$$

Write $Q_s$ for the projection onto $H(s)$, such that $Q=\sum_sQ_s$. Then
$$Q_s(\Phi_j(Y)\otimes
I_H)|_{E^{\otimes s}\otimes_\pi H}=Q_s(T_{\xi^{(j)}_0}\otimes I_H)|_{E^{\otimes s}\otimes_\pi H}.$$
\end{prop}

\noindent{Proof.} Let $\Phi_0(Y)\otimes I_H=diag(\pi(a_0),\phi_1(a_1)\otimes I_H,...)$.
For every $s\geq 0$ the subspace $H(s)$ is $\Phi_0(Y)\otimes I_H$-invariant and so is $\mathfrak{M}$. Since $(\Phi_0(Y)^*\otimes I_H)|_{\mathfrak{M}}=(\phi_{\infty}(a_0^*)\otimes I_H)|_{\mathfrak{M}}$ we get
$$(\Phi_0(Y)\otimes I_H)|_{\mathfrak{M}}=(\phi_{\infty}(a_0)\otimes I_H)|_{\mathfrak{M}}.$$

2) Set $S^*\otimes I_H=\Phi_j(Y)^*\otimes I_H-T^*_{\xi^{j}_0}\otimes I_H$. Then $(S^*\otimes I_H)|_{\mathfrak{M}}=0$, that is $Q(S^*\otimes I_H)Q=0$. Hence, $Q(S\otimes I_H)Q=0$ and, in particular, $Q_{s+j}(S\otimes I_H)Q_s=0$ for every $s\geq 0$.

\begin{flushright}
$\square$
\end{flushright}

\begin{cor}\label{a_0=a_1}  Let $Z_0\otimes I_H$ be an operator in $Alg\ Lat\
\rho_\pi(H^{\infty}(E))$, admitting a representation
$$Z_0\otimes I_H=diag(\pi(a_0),\phi_1(a_1)\otimes I_H,...),\ \ a_s\in M.$$
Then
$$(Z_0\otimes I_H)|_{\mathfrak{M}}=(\phi_{\infty}(a_0)\otimes
I_H)|_{\mathfrak{M}}.$$
Hence, for all  $s\geq 0$
$$(Z_0\otimes I_H)|_{H(s)}=(\phi_{\infty}(a_0)\otimes
I_H)|_{H(s)},$$
and
$$a_0=a_1.$$
\end{cor}
\noindent{Proof.} From the previous proposition we get $(Z_0\otimes I_H)|_{\mathfrak{M}}=(\phi_{\infty}(a_0)\otimes I_H)|_{\mathfrak{M}}$ and, for every $s\geq 0$, $(Z_0\otimes I_H)|_{H(s)}=(\phi_{\infty}(a_0)\otimes I_H)|_{H(s)}$.
As we saw $H(0)=H$ and $H(1)=E\otimes_\pi H$, and the equality $a_0=a_1$ follows.
\begin{flushright}
$\square$
\end{flushright}

At this point we are able to show that for every $Y\otimes I_H\in Alg\ Lat\ \rho_\pi(H^{\infty}(E))$, its ``0"-th Fourier coefficient $\Phi_0(Y)\otimes I_H$ is in $\rho_\pi(H^{\infty}(E))$ (Theorem $\ref{Phi_0 in H infty}$). We will use the following two simple observations.
Let $\mathcal{A}$ be an operator algebra acting on Hilbert space $H$, $\mathcal{N}\in Lat\ \mathcal{A}$ and $T\in Alg\ Lat\ \mathcal{A}$. Then $T|_{\mathcal{N}}\in Alg\ Lat\ (\mathcal{A}|_{\mathcal{N}})$. If $\mathcal{B}$ another operator algebra acting on a Hilbert space $K$ such that there is a unitary $W:H\rightarrow K$ with $W\mathcal{A}W^*=\mathcal{B}$, then for every $T\in Alg\ Lat\ \mathcal{A}$, $WTW^*\in Alg\ Lat\ \mathcal{B}$.

The first claim is evident. For the second note that if
$A=W^*BW$ for some $B\in\mathcal{B}$ and if $\mathcal{M}\in Lat\ \mathcal{B}$, then
$$A(W^*\mathcal{M})=W^*BW(W^*\mathcal{M})=W^*B\mathcal{M}\subset W^*\mathcal{M}\subset \mathcal{M}^*.$$
Thus, $W^*\mathcal{M}\in Lat\ \mathcal{A}$.
Now if $T\in Alg\ Lat\ \mathcal{A}$ then $TW^*\mathcal{M}\subset W^*\mathcal{M}$, hence $(WTW^*)\mathcal{M}\subset \mathcal{M}$.

Fix $s\geq 1$ and set $K_s:=E^{\otimes s}\otimes_\pi H$ and $\pi_s:=\phi_s\otimes I_H$. Then $\pi_s$ is a faithful normal representation of $M$ on $K_s$. Then the space $\mathcal{F}(E)\otimes_{\pi_s}K_s$ may be identified with the subspace $G_s:=\sum^{\oplus}_{l\geq s}E^{\otimes l}\otimes_{\pi}H\subset \mathcal{F}(E)\otimes_{\pi}H$ as follows.

Let $k_s=\zeta\otimes h\in K_s$, with $\zeta\in E^{\otimes s}$. Then the formula
$$W_{s,l}:\xi_1\otimes...\otimes\xi_l\otimes\zeta\otimes h\mapsto \xi_1\otimes...\otimes\xi_l\otimes k_s$$
defines the identification
$$ E^{\otimes l+s}\otimes_\pi
H\cong E^{\otimes l}\otimes_{\pi_s}K_s,\,\,\, l,s\geq 0.$$
Hence, we obtain a unitary operator
$$W_s:G_s=\sum^{\oplus}_{l\geq s}E^{\otimes l}\otimes_{\pi}H\rightarrow
\mathcal{F}(E)\otimes_{\pi_s}K_s.$$

Write $\rho_0$ for the induced representation $\pi^{\mathcal{F}(E)}$ and $\rho_s$ for the induced representation $\pi_s^{\mathcal{F}(E)}$ of $H^{\infty}(E)$ on $\mathcal{F}(E)\otimes_{\pi_s}K_s$. Thus, $\rho_s(X)=X\otimes I_{K_s}$ when $X\in H^{\infty}(E)$.

Further, by $\rho_0|_{G_s}$ we denote the representation obtained by the restriction of $\rho_0$ to $G_s$. Thus, $\rho_0|_{G_s}(X)=(X\otimes I_H)|_{G_s}$ when $X\in H^{\infty}(E)$.

\begin{lem}\label{Identification rho_0|G_s with rho_s} Fix  $s\geq 0$ and let $W_s$ be as above. Then

1) $W_s^*\rho_s(X)W_s=\rho_0(X)|_{G_s}$, $X\in H^{\infty}(E)$;

2) for every $Y\otimes I_H\in Alg\ Lat\ \rho_0(H^{\infty}(E))$, the restriction $(Y\otimes I_H)|_{G_s}$ is in $Alg\ Lat\ (\rho_0(H^{\infty}(E))|_{G_s})$, and $W(Y\otimes I_H)|_{G_s}W^*\in Alg\ Lat\ \rho_s(H^{\infty}(E))$.
\end{lem}
\noindent{Proof.} 1) Enough to check it for the generators of the Hardy algebra. So, take $\phi_{\infty}(a)\in H^{\infty}(E)$, $a\in M$. For $\mathfrak{z}\otimes\zeta\otimes h\in E^{\otimes l+s}\otimes_\pi H$, where we write $\mathfrak{z}$ for $\xi_1\otimes...\otimes\xi_l$ and $\zeta\in E^{\otimes s}$, we have
\begin{align*}
W_s^*\rho_0(\phi_{\infty}(a))W_s(\mathfrak{z}\otimes\zeta\otimes h)=W_s^*(\phi_{\infty}(a)\otimes I_{K_s})(\mathfrak{z}\otimes k_s)=\\
(\phi_l(a)\mathfrak{z})\otimes\zeta\otimes h=(\phi_{\infty}(a)\otimes I_H)(\mathfrak{z}\otimes\zeta\otimes h)=\\
\rho_0(\phi_{\infty}(a))(\mathfrak{z}\otimes\zeta\otimes h).
\end{align*}
Now take $T_\xi\in H^{\infty}(E)$, $\xi\in E$. Then
\begin{align*}W_s^*\rho_s(T_{\xi})W_s(\mathfrak{z}\otimes\zeta\otimes h)=W_s^*(T_{\xi}\otimes I_{K_s})(\mathfrak{z}\otimes k_s)=\\
\xi\otimes \mathfrak{z}\otimes\zeta\otimes h=
(T_{\xi}\otimes I_{H})(\mathfrak{z}\otimes\zeta\otimes h)=\\
\rho_0(T_\xi)(\mathfrak{z}\otimes\zeta\otimes h),
\end{align*}

Hence, $$W_s^*\rho_s(X)W_s=\rho_0(X)|_{G_s},\,\,\, X\in H^{\infty}(E).$$

2) Let $Y\otimes
I_H\in Alg\ Lat\ \rho_0(H^{\infty}(E))$. Then using by the observations done before lemma, $(Y\otimes I_H)|_{G_s}\in Alg\ Lat\ (\rho_0(H^{\infty}(E))|_{G_s})$ and $W(Y\otimes I_H)|_{G_s}W^*\in Alg\ Lat\ \rho_s(H^{\infty}(E))$.

\begin{flushright}
$\square$
\end{flushright}

\begin{thm}\label{Phi_0 in H infty}
Let $Y\otimes I_H\in Alg\ Lat\ \rho_0(H^{\infty}(E))$. Then
$$\Phi_0(Y)\otimes I_H\in \rho_0(H^{\infty}(E)).$$

\end{thm}
\noindent{Proof.} Write $Z_0\otimes I_H$ for $\Phi_0(Y)\otimes I_H$. $Z_0\otimes I_H$ admits the representation
$$Z_0\otimes I_H=diag(\pi(a_0),\pi_1(a_1),...),\ \ a_i\in M.$$
Thus, the operator $Z_0\otimes I_H$ is as in Corollary $\ref{a_0=a_1}$, 1). Hence
$$a_0=a_1.$$
Fix now $s\geq 1$ arbitrary. Then
$$(Z_0\otimes I_H)|_{G_s}=diag(\pi_s(a_s),\pi_{s+1}(a_{s+1}),...).$$
By Lemma $\ref{Identification rho_0|G_s with rho_s}$,
$W_s(Z_0\otimes I_H)|_{G_s}W_s^*\in Alg\ Lat\ \rho_s(H^{\infty}(E))$, and has the matrix representation
$$W_s(Z_0\otimes I_H)|_{G_s}W_s^*=diag(\phi_s(a_s)\otimes
I_{K_s},\phi_1(a_{s+1})\otimes I_{K_s},...).$$
The operator $W_s(Z_0\otimes I_H)|_{G_s}W_s^*$ satisfies all conditions of Corollary $\ref{a_0=a_1}$, where we replace $\pi$ by $\phi_s\otimes I_{K_s}$. Hence $a_s=a_{s+1}$.

Since the choice of $s$ is arbitrary, we obtain  $a_0=a_1=a_2=...$.
Thus, $Z_0\otimes I_H=\Phi_j(Y)\otimes I_H=\phi_{\infty}(a_0)\otimes I_H$.
\begin{flushright}
$\square$
\end{flushright}


\begin{rem}\label{Action of Z_j on big M_s}
For the Fourier coefficient $\Phi_j(Y)\otimes I_H$ with $j\geq 1$, the above technique is not working and we obtain only the following. Write $Z_j\otimes I_H$ for the $j$-th Fourier coefficient $\Phi_j(Y)\otimes I_H$, $j\geq 1$, of some $Y\otimes I_H\in Alg\ Lat\ \rho_0(H^{\infty}(E))$. According to Theorem
$\ref{Repres. of Phi_j}$, $Z_j\otimes I_H$ represented as a multiplication by sequence $\{\xi^{(j)}_0,\xi^{(j)}_1,\xi^{(j)}_2,...\}\subset E^{\otimes j}$. Fix some $s\geq 0$. Then, as it is easy to see from the definition of $G_s$, the restriction $(Z_j\otimes I_H)|_{G_s}$  is represented by multiplication by sequence $\{\xi^{(j)}_s,\xi^{(j)}_{s+1},...\}\subset E^{\otimes j}$. Let further, $\mathfrak{M}_s$ be the symmetric $\rho_s(H^{\infty}(E))$-coinvariant subspace of $\mathcal{F}(E)\otimes_{\pi_s}K_s$, and write $R_s$ for the projection onto $\mathfrak{M}_s$ (with $R_0=Q$). Then, by Lemma $\ref{Identification rho_0|G_s with rho_s}$, and by Theorem $\ref{Action of Phi_j on big M}$, applied to the representation $\rho_s$ and the subspace $\mathfrak{M}_s$, we obtain that
\begin{equation}\label{Compression Phi_j on big M_s}
R_sW_s(Z_j\otimes I_H)W_s^*|_{\mathfrak{M}_s}=R_sW_s(T_{\xi^{(j)}_s}\otimes I_H)W_s^*|_{\mathfrak{M}_s}.
\end{equation}
\end{rem}
~

In the following we shall drop the upper index for a general element $\xi\in E^{\otimes j}$.

Now we may prove the reflexivity of $\rho_\pi(H^{\infty}(E))$ where $\pi: M\rightarrow B(H)$ is a reducible representation of a factor. In the proof we shall use the following lemma.
\begin{lem}\label{pi(a)h'=0 impl a=0} Let $M$ be a factor and $\pi:M\rightarrow B(H)$ be a faithful normal representation of $M$ on $H$. Let $H'\subseteq H$ be a nontrivial $\pi(M)$-invariant subspace. Then

1) If $a\in M$ and $\pi(a)|_{H'}=0$ then $a=0$;

2) Let $k\geq 0$, $j\geq 1$ and $\xi\in E^{\otimes j}$ such that $\xi\otimes E^{\otimes k}\otimes_\pi H'=\{0\}$. Then $\xi=0$.
\end{lem}
\noindent{Proof.} 1) Let $a\in M$ and $\pi(a)h'=0$ for every $h'\in H'$. Write $P'$ for the projection onto $H'$. Clearly, $P'\in \pi(M)'$. For every $l\in H$ and $h'\in H'$ one have $\langle l,\pi(a)h'\rangle=0$. Hence, $\pi(a)P'=0$. But the same is true if we replace $P'$ by its central carrier $C_{P'}$. Since $M$ is factor, so is $\pi(M)'$ and $C_{P'}=I_H$. Thus, $\pi(a)=0$, and since $\pi$ is faithful we get $a=0$.

2) We distinguish two cases $k=0$ and $k\geq 1$.

Let $k=0$. If $\xi\otimes_\pi H'=\{0\}$ then for $\zeta\in E^{\otimes j}$, $h'\in H$ and $l\in H$,
$$\langle \zeta\otimes l,\xi\otimes h'\rangle=0.$$
Then
$$\langle l,\pi(\langle \zeta,\xi\rangle)h'\rangle=0,$$

Since $l\in H$ is arbitrary, we obtain
$$\pi(\langle\zeta,\xi\rangle)h'=0.$$
Write $Q'_0$ for the projection in $H$ onto $H'$. Clearly, $Q'_0\in \pi(M)'$.
We get
$$\pi(\langle\zeta,\xi\rangle))Q'_0=0.$$
As in 1), since $M$ is factor so is $\pi(M)'$ and the central carrier $C_{Q'_0}=I_H$. Hence, $\pi(\langle\zeta,\xi\rangle)=0$.
But $\pi$ is faithful and $\zeta\in E^{\otimes j}$ is arbitrary. Thus $\xi=0$.

Let now $k\geq 1$. Then $\xi\otimes E^{\otimes k}\otimes_\pi H'=\{0\}$. For $\zeta\in E^{\otimes j}$, $\theta_1,\theta_2\in E^{\otimes k}$, $h'\in H$ and $l\in H$ we have

$$\langle \zeta\otimes\theta_2\otimes l,\xi\otimes \theta_1\otimes h'\rangle=0.$$
Hence,
$$\langle\theta_2\otimes l,(\phi_k(\langle \zeta,\xi\rangle)\theta_1)\otimes h'\rangle=0.$$
Since $\theta_2\in E^{\otimes k}$ and $l\in H$ are arbitrary, we obtain
$$(\phi_k(\langle\zeta,\xi\rangle)\theta_1)\otimes h'=0.$$

Write $Q'$ for the projection in $E^{\otimes k}\otimes_\pi H$ onto $E^{\otimes k}\otimes_\pi H'$. Clearly, $Q'\in (\phi_k(M)\otimes I_H)'$ and
$$(\phi_k(\langle\zeta,\xi\rangle)\otimes I_H)Q'=0.$$

As in case when $k=0$, since $M$ is factor so is $(\phi_k(M)\otimes I_H)'$, and the central carrier $C_{Q'}=I_{E^{\otimes k}\otimes_\pi H}$. Hence, $\phi_k(\langle\zeta,\xi\rangle)\otimes I_H=0$.
But $\phi_k$ is faithful and $\zeta\in E^{\otimes j}$ is arbitrary. Thus $\xi=0$.
\begin{flushright}
$\square$
\end{flushright}

\begin{thm} \label{Reflexivity of Hardy alg over factor.}
Let $M$ be any factor and let $\pi:M\rightarrow B(H)$ be a reducible faithful normal representation of $M$ on a Hilbert space $H$. Then the algebra $\rho_\pi(H^{\infty}(E))$ is reflexive.
\end{thm}
\noindent{Proof.} Let $Y\otimes I_H\in Alg\ Lat\ \rho_\pi(H^{\infty}(E))$. We already have seen (even in a more general situation) that $\Phi_0(Y)\otimes I_H \in\rho_\pi(H^{\infty}(E))$. Thus we need to show that, under our assumptions, $\Phi_j(Y)\otimes I_H \in \rho_\pi(H^{\infty}(E))$, for every $j\geq 1$.

Let $M$ be an arbitrary factor. Since $\pi$ is reducible, the Hilbert space $H$ splits onto direct sum $H=H_0\oplus H_1$ of $\pi(M)$-invariant subspaces $H_0, H_1\neq(0)$, each of them is wandering with respect to the covariant representation $(V,\sigma)$, $\sigma(\cdot)=\phi_\infty(\cdot)\otimes I_H$.

Set $\mathcal{L}_k:=\mathfrak{L}^k(H_1)=E^{\otimes k}\otimes_\pi H_1\subset \mathcal{F}(E)\otimes_\pi H$, where $\mathfrak{L}$ is the generalized shift associated with $(V,\sigma)$. Since $(V,\sigma)$ is isometric, it is clear that $E^{\otimes k}\otimes_\pi H_1\perp E^{\otimes k}\otimes_\pi H_0$.
Now set $\mathcal{N}_k:=H_0\oplus \mathcal{L}_k$. Since $H_0$ and $\mathcal{L}_k$ both are wandering and $H_0\perp H_1$, $\mathcal{N}_k$ is a wandering subspace. $\mathcal{N}_k$ generates the subspace in $\mathcal{F}(E)\otimes_\pi H$
$$\mathcal{M}_{\mathcal{N}_k}:=\sum_s^{\oplus}\mathfrak{L}^s(\mathcal{N}_k),$$
which is unitarily isomorphic to the space
$$\mathcal{F}(E)\otimes_{\sigma'}\mathcal{N}_k=\mathcal{N}_k\oplus E\otimes_{\sigma'}\mathcal{N}_k\oplus...,$$
where $\sigma':=\sigma|_{\mathcal{N}_k}$. It should be noted that in our assumption the representation $\sigma'$ is faithful. This follows from Lemma $\ref{pi(a)h'=0 impl a=0}$, 1), since $\sigma'|_{H}=\pi|_H$ and if $\pi(a)\mathcal{N}_k=0$ for some $a\in M$, then $\pi(a)|_{H_0}=0$ for the nontrivial $\pi(M)$-invariant subspace $H_0\subset H$. By Lemma $\ref{pi(a)h'=0 impl a=0}$, 1), $a=0$.

Clearly, $\mathcal{M}_{\mathcal{N}_k}$ is a $\rho_\pi(H^{\infty}(E))$-invariant subspace, and we denote by $V_k$ the unitary from $\mathcal{M}_{\mathcal{N}_k}$ onto $\mathcal{F}(E)\otimes_{\sigma'}\mathcal{N}_k$.

Write $Z_j\otimes I_H$ for the operator $\Phi_j(Y)\otimes I_H$, $Y\otimes I_H\in Alg\ Lat\ \rho_\pi(H^{\infty}(E))$. Then the restriction $(Z_j\otimes I_H)|_{\mathcal{M}_{\mathcal{N}_k}}$ is in $Alg\ Lat\ \rho_\pi(H^{\infty}(E))|_{\mathcal{M}_{\mathcal{N}_k}}$ and the operator $V_k(Z_j\otimes I_H)V_k^*$ is in $Alg\ Lat\ \rho_{\sigma'}(H^{\infty}(E))$. Applying Corollary $\ref{Alg Lat in commutant of M'}$ to $V_k(Z_j\otimes I_H)V_k^*$, we deduce that there is $Z_j'\in \mathcal{L}(\mathcal{F}(E))$ such that
$$Z_j'\otimes I_{\mathcal{N}_k}=V_k(Z_j\otimes I_H)V_k^*\in Alg\ Lat\ \rho_{\sigma'}(H^{\infty}(E)).$$
By Theorem $\ref{Repres. of Phi_j}$ on the representation of the Fourier coefficients we obtain that there are $\mathfrak{z}_{i}\in E^{\otimes j}$, $i=0, 1,2,...$, such that
$$(\Phi_j(Z_j')\otimes I_{\mathcal{N}_k})|_{E^{\otimes k}\otimes_{\sigma'}\mathcal{N}_k}=(T_{\mathfrak{z}_{k}}\otimes I_{\mathcal{N}_k})|_{E^{\otimes k}\otimes_{\sigma'}\mathcal{N}_k}, \,\,\, k\geq 0.$$
Recall the formula $\Phi_j(Z_j')\otimes I_{\mathcal{N}_k}=(\sum_sP_{s+j}Z_j'P_s)\otimes I_{\mathcal{N}_k}$, where, as usual, $P_s$ is the projection $P_{E^{\otimes s}}$, and write $R_s=P_s\otimes I_{\mathcal{N}_k}$.

Take $f_0\in H_0$ and $\theta\otimes l\in \mathcal{L}_k$ arbitrary.
Then
$$(\Phi_j(Z_j')\otimes I_{\mathcal{N}_k})(f_0+\theta\otimes l)=\mathfrak{z}_{0}\otimes(f_0+\theta\otimes l).$$
On the other hand
$$(\Phi_j(Z_j')\otimes I_{\mathcal{N}_k})(f_0+\theta\otimes l)=R_{j}(V_{k}(Z_j\otimes I_H)V_k^*)R_0(f_0+\theta\otimes l)=\xi_0\otimes f_0+\xi_k\otimes\theta\otimes l,$$
where $\{\xi_k\}_{k\geq 0}\subset E^{\otimes j}$ is the sequence which corresponds to $Z_j\otimes I_H$ by Theorem $\ref{Repres. of Phi_j}$.
Thus, $$(\Phi_j(Z_j')\otimes I_{\mathcal{N}_k})(f_0+\theta\otimes l)=(Z_j\otimes I_H)(f_0+\theta\otimes l)$$
or
\begin{equation}\label{L_zita_0=Phi_j on N_j}
\mathfrak{z}_{0}\otimes(f_0+\theta\otimes l)=\xi_0\otimes f_0+\xi_k\otimes \theta\otimes l.
\end{equation}
Set $\zeta_0:=\mathfrak{z}_{0}-\xi_0$. Hence
$$L_{\zeta_0}|_{H_0}=0.$$
Using Lemma $\ref{pi(a)h'=0 impl a=0}$, 2), taking there $k=0$, we obtain $\zeta_0=0$. Thus, $\xi_0=\mathfrak{z}_{0}$.

From $(\ref{L_zita_0=Phi_j on N_j})$ we obtain now
$$L_{\xi_k}|_{\mathcal{L}_k}=L_{\xi_0}|_{\mathcal{L}_k}.$$
Set $\zeta_1:=\xi_k-\xi_0$. Thus $L_{\zeta_1}|_{\mathcal{L}_k}=0$.
Again, using Lemma $\ref{pi(a)h'=0 impl a=0}$, 2) (with $k\geq 1$), we get
$\zeta_1=0$, i.e. $\xi_k=\xi_0$.
Since $k\geq 1$ is arbitrary, we obtain that $\xi_0=\xi_1=...$.

This shows that
$\Phi_j(Y)\otimes I_H\in\rho_\pi(H^{\infty}(E))$ for any $j\geq 0$, hence, $Y\otimes I_H\in \rho_\pi(H^{\infty}(E))$, i.e. $\rho_\pi(H^{\infty}(E))=Alg\ Lat\ \rho_\pi(H^{\infty}(E))$.
\begin{flushright}
$\square$
\end{flushright}

\begin{cor}\label{Reflexivity over factors II,III} If $M$ is a factor of type $II$ or $III$, then for every faithful normal representation $\pi: M\rightarrow B(H)$, the algebra $\rho_\pi(H^{\infty}(E))$ is reflexive.
\end{cor}
\noindent{Proof.} If $\pi$ is irreducible, then $\pi(M)'=\mathbb{C}I$. Hence, $\pi(M)=\pi(M)''=B(H)$, i.e. $\pi(M)$, and therefore $M$, is of type $I$. When this is not the case, $\pi$ is reducible and the previous theorem applies.
\begin{flushright}
$\square$
\end{flushright}

\begin{rem} 1) In ~\cite{APo} Arias and Popescu proved the reflexivity of $\rho(H^{\infty}(\mathbb{C}^n))$ over the factor (of type $I$) $M=\mathbb{C}$, without the assumption of reducibility of $\pi$.

2) In the previous chapter we saw that if $M$ is a factor of type $III$ then $\rho_\pi(H^{\infty}(E))$ is even hyperreflexive.
\end{rem}
\begin{flushright}
$\square$
\end{flushright}

\subsection{Analytic crossed product. An example.}

Here we consider the $W^*$-correspondence  $_{\alpha}M$ over a $W^*$-algebra $M$. We recall that
$\alpha\in End(M)$ and that $p=\alpha(1)$ is a projection in $M$. We set $E=M\otimes_{\alpha}M$, that is the self-dual completion of the algebraic tensor product $M\otimes M$ of algebras with the relations $ac\otimes b=a\otimes\alpha(c)b$ and equipped with the inner product defined by $\langle a\otimes b,c\otimes d\rangle=b^*\alpha(a^*c)d$.
We identify $M\otimes_{\alpha}M$ with $\overline{\alpha(M)M}$ and it is easy to see that the inner product in $\overline{\alpha(M)M}$ turns to be $\langle\alpha(a)b,\alpha(c)d\rangle=(\alpha(a)b)^*\alpha(c)d$.
The left action of $M$ on $E$ we define by $a\cdot \xi=\alpha(a)\xi,$
where $a\in M$ and $\xi\in E$. So, $E$ has the structure of a $W^*$-correspondence over $M$ and we write $_{\alpha}M$ for it.
In what follows we always assume $\alpha$ to be injective.
The map $a\mapsto \alpha(1)a=pa$ gives us the identification of $_{\alpha}M$ with $\alpha(1)M$. Thus, $E=\ _{\alpha}M=\alpha(1)M=pM$.

Note that for every $k\geq 0$ one have $\alpha^k(1)=\alpha^{k-1}(p)\alpha^{k-2}(p)...\alpha(p)p$, and $\alpha^k(p)\geq\alpha^{k+1}(p)$, for every $k\geq 0$.

For every $k\geq 0$, the map
$\xi_1\otimes...\otimes \xi_k\mapsto \alpha^{k-1}(\xi_1)...\alpha(\xi_{k-1})\xi_k,$
gives an identification $E^{\otimes k}=\ _{\alpha^k}M$, $_{\alpha^k}M=\alpha^k(1)M$. The left action of $M$ on $_{\alpha^k}M$ is given by the formula $a\cdot\xi_k=\alpha^k(a)\xi_k$, $a\in M$ and $\xi_k\in\ _{\alpha^k}M$
The full Fock space over $_{\alpha}M$ is now

$$\mathcal{F}(E)=\mathcal{F}( _{\alpha}M)=\sum^{\oplus}_{k\geq 0}\, _{\alpha^k}M.$$

For the left action $\phi_\infty$ of $M$ on $\mathcal{F}(_{\alpha}M)$ we write $\alpha_\infty$, which is now given by
$\alpha_{\infty}=diag(\alpha^{0},\alpha,\alpha^{2},...)$. Thus, if $(\xi_k)\in\mathcal{F}(_\alpha M)$ then $\phi_\infty(a)(\xi_k)=\alpha_\infty(a)(\xi_k)=(\alpha^{k}(a)\xi_k)$.

Let $\xi\in E=\ _{\alpha}M$ and $\xi_k\in\ _{\alpha^k}M$, then
$T_{\xi}\xi_k=\xi\otimes \xi_k=\alpha^{k}(\xi)\xi_{k}\in\ _{\alpha^{k+1}}M$. Clearly,
$T_{\alpha(a)\xi b}=\alpha_{\infty}(a)T_{\xi}\alpha_{\infty}(b)$. Since every $\xi\in E$ has a form
$\xi=\alpha(1)a$, $a\in M$, we get $T_{\xi}=T_{\alpha(1)}\alpha_{\infty}(a)$. Thus, the operator $T_{\xi}$ is completely determined by $T_{\alpha(1)}=T_p$. Every $\xi_k\in\ _{\alpha^k}M$ has the form $\alpha^k(1)a_k$, $a_k\in M$. Then
$T_{\xi}\xi_k=\alpha^{k}(\xi)\xi_{k}=
\alpha^{k+1}(1)\alpha^k(a)a_k$.
So, the set of generators of $H^{\infty}( _{\alpha}M)$ consists of left multiplications $\alpha_{\infty}(a)$ by elements $a$ of $M$, and the creation operator $T_{\alpha(1)}$.

Let $\pi:M\rightarrow B(H)$ be a faithful normal representation
of $M$ on Hilbert space $H$.
For every $k\geq 1$, the space $_{\alpha^k}M\otimes_{\pi}H$ can be identified with $\pi(\alpha^k(1))H=\pi(\alpha^{k-1}(p))H$ via
$$\xi_1\otimes...\otimes\xi_k\otimes h\mapsto\pi(\alpha^{k-1}(\xi_1)...\alpha(\xi_{k-1})\xi_k)h.$$
For $k=0$ this formula reduces to $M\otimes_{\pi}H\cong\overline{\pi(M)H}$.
So, the space $\mathcal{F}(_{\alpha}M)\otimes_{\pi}H$ can be identified with a subspace of $l^{2}(\mathbb{Z}_+, H)$:
$$\mathcal{F}(_{\alpha}M)\otimes_{\pi}H\cong \sum^{\oplus}_{k\geq 0}\pi(\alpha^{k}(1))H.$$
Let us consider the induced representation $\rho(H^{\infty}(_{\alpha}M))=
\pi^{\mathcal{F}(_{\alpha}M)}(H^{\infty}(_{\alpha}M))$. If $Y\otimes I_H\in Alg\ Lat(\rho(H^{\infty}(_{\alpha}M)))$ then, as it follows from Theorem $\ref{Repres. of Phi_j}$, the $j$-th Fourier coefficient is represented as multiplication by  the matrix
(we omit the upper indices in $L$ and the upper indices in $\xi\in E^{\otimes j}$)
$$\left(
  \begin{array}{cccccc}
    0   & 0   & 0   & 0   & ... & ... \\
    ... & ... & ... & ... & ... & ... \\
    0   & ... & ... & ... & ... & ... \\
    L_{\xi_0}   & 0   & 0   & 0 & ... & ... \\
    0   & L_{\xi_1}   & 0   & 0 & ... & ... \\
    0   & 0   & L_{\xi_2}   & 0 & ... & ... \\
    ... & ... & ... & ... & ... & ... \\
  \end{array}
\right),$$
where each $\xi_s\in _{\alpha^j}M$. Hence there exist sequence $\{a_s\}_{s\geq 0}\subset M$ such that $\xi_s=\alpha^j(1)a_s$.
Notice also that the adjoint operator $(T_\xi\otimes I_H)^*$ acts on $_{\alpha^k}M\otimes_\pi H$ by the formula
$$T_\xi^*\otimes I_H:\xi_k\otimes h=\pi(\alpha^k(1)a_k)h\mapsto \pi(\alpha^{k-1}(a^*))\pi(\alpha^{k}(1)a_{k})h,$$
$\xi=\alpha(1)a$, $\xi_{k}=\alpha^{k}(1)a$.

By Theorem $\ref{Reflexivity of Hardy alg over factor.}$, if $\pi: M\rightarrow B(H)$ is a reducible representation of a factor (in particular if $M$ is a factor of type $II$ or $III$), then the algebra $\rho_\pi(H^{\infty}(_{\alpha}M))$ is reflexive, for every endomorphism $\alpha\in End(M)$.
In ~\cite{Ka} E. Kakariadis, showed that in the special case when $\alpha$ is a unitarily implemented automorphism of $M$, the algebra $\rho(H^{\infty}( _{\alpha}M))$ is reflexive. We want first to show how Theorem $\ref{Repres. of Phi_j}$ can be used to give another proof of the reflexivity of the analytic crossed product in this case.

So, let $v\in U(M)$ be some unitary in $U(M)$, the unitary group of $M$. Let $\alpha\in Aut (M)$, and assume that $\alpha(a)=vav^*$.
Write $u=\pi(v)\in U(H)$, where $U(H)$ is the unitary group of the Hilbert space $H$. Hence, $\pi(\alpha(a))=u\pi(a)u^*$.

Fix some $h\in H$, some $0< r <1$ and set
$$\mathcal{M}:=close\{\pi(b)h+\pi(vb)rh+\pi(v^2b)r^2h+...:b\in M\}.$$
It is easy to see that $\mathcal{M}$ is a $\rho(H^{\infty}( _{\alpha}M))$-coinvariant subspace. Indeed,
let $x=\pi(b)h+\pi(vb)rh+\pi(v^2b)r^2h+...\in \mathcal{M}$ and let $\xi=a\in\ _{\alpha}M$, then
$$(T^*_{\xi}\otimes I_H)x=\pi(a^*)\pi(vb)rh+\pi(\alpha(a^*))\pi(v^2b)r^2h+...+\pi(\alpha^{s-1}(a^*))\pi(v^s)r^sh+...,$$
and since $\pi(\alpha^{s-1}(a^*))\pi(v^s)r^sh=\pi(v^{s-1}a^*(v^*)^{ s-1}v^s)r^sh=\pi(v^{s-1}a^*vr)r^{s-1}h$, we obtain that for $c=a^*vr$
$$(T^*_{\xi}\otimes I_H)x=\pi(c)h+\pi(vc)rh+\pi(v^2c)r^2h+...\in \mathcal{M}.$$
If $a\in M$ then $(\alpha_s(a^*)\otimes I_H)\pi(v^sb)r^sh=\pi(\alpha^s(a^*)v^sb)r^sh=\pi(v^sa^*v^{*s}v^sb)r^sh=\pi(v^sa^*b)r^sh$. Hence,
$$(\alpha_{\infty}(a^*)\otimes I_H)x=\pi(a^*b)h+\pi(va^*b)rh+\pi(v^2a^*b)r^2h+....\in \mathcal{M}.$$
Thus $\mathcal{M}$ is indeed coinvariant.

Take $$y:=h+\pi(v)rh+\pi(v^2)r^2h+...+\pi(v^s)r^sh+...\in \mathcal{M}.$$
Then, using the notation of Theorem $\ref{Repres. of Phi_j}$,
$$(\Phi_j(Y)\otimes I_H)^*y=L^*_{\xi_0}\pi(v^j)r^jh+L^*_{\xi_1}\pi(v^{j+1})r^{j+1}h+...+
L^*_{\xi_s}\pi(v^{j+s})r^{j+s}h+...,$$
where $\xi_s\in E^{\otimes j}$ and $L_{\xi_s}=(T_{\xi_s}\otimes I_H)|_{E^{\otimes s}\otimes_{\pi}H}:E^{\otimes s}\otimes_{\pi}H\rightarrow E^{\otimes s+j}\otimes_{\pi}H$. Since $(\Phi_j(Y)\otimes I_H)^*y\in \mathcal{M}$, there exist some sequence $\{c_i\}_i\subset M$ such that
$$(\Phi_j(Y)\otimes I_H)^*y=\lim_i[h+\pi(c_i)h+\pi(vc_i)rh+\pi(v^2c_i)r^2h+...+\pi(v^sc_i)r^sh+...].$$

Hence,
$$\lim_i\pi(c_i)h=L^*_{\xi_0}\pi(v^j)r^jh=\pi(a_0^*v^j)r^jh,$$
and
$$\lim_i\pi(v^sc_i)r^sh=L^*_{\xi_s}\pi(v^{j+s})r^{j+s}h=\pi(\alpha^s(a_s^*)v^{j+s})r^{j+s}h,$$

Then,
$$\lim_i\pi(c_i)h=\pi(a_0^*v^j)r^jh=\pi(v^j)\pi(\alpha^{-j}(a_0^*))r^jh,$$
and
$$\lim_i\pi(v^sc_i)r^sh=\pi(\alpha^s(a_s^*)v^{j+s})r^{j+s}h=
\pi(v^{j+s})\pi(\alpha^{-(j+s)}(\alpha^s(a_s^*)))r^{j+s}h.$$
But $\lim_i\pi(v^sc_i)r^sh=\pi(v^s)r^s\lim_i\pi(c_i)h=
\pi(v^{j+s})\pi(\alpha^{-j}(a_0^*))r^{j+s}h$, and we get

$$\pi(v^{j+s})\pi(\alpha^{-j}(a_0^*))r^{j+s}h=\pi(v^{j+s})\pi(\alpha^{-(j+s)}(\alpha^s(a_s^*))r^{j+s}h.$$

Since $\pi(v^k)$ is unitary for every $k$, $r$ is scalar and this equality holds for every $h\in H$, we obtain

$$\pi(\alpha^{-j}(a_0^*))=\pi(\alpha^{-(j+s)}(\alpha^s(a_s^*)).$$
Since $\pi$ is faithful it follows that
$$\alpha^{-j}(a_0^*)=\alpha^{-(j+s)}(\alpha^s(a_s^*))=\alpha^{-j}(a_s^*).$$
Thus, $a_s=a_0$ for every $s\geq 0$.
We showed that if $Y\otimes I_H\in Alg\ Lat\ \rho(H^{\infty}( _{\alpha}M))$ then for every $j\geq 1$, $\Phi_j(Y)\otimes I_H$ is in $\rho(H^{\infty}( _{\alpha}M))$. Hence we proved
\begin{thm} Let $\pi:M\rightarrow B(H)$ be a faithful normal representation of the $W^*$-algebra $M$. Then the algebra $\rho_\pi(H^{\infty}( _{\alpha}M))$ is reflexive, whenever $\alpha$ is an automorphism that is unitarily implemented.
\end{thm}
\begin{flushright}
$\square$
\end{flushright}

\begin{cor}\label{ReflexivityOfAnalyt.CrossedProdOverFactorsForEveryAutom} If $M$ is a factor and $\alpha\in Aut(M)$, then $\rho_\pi(H^{\infty}(\ _{\alpha}M))$ is reflexive for any normal representation $\pi$.
\end{cor}
\noindent{Proof.}  The type $II$ and type $III$ cases follow from Corollary $\ref{Reflexivity over factors II,III}$.
If $M$ is a type $I$ factor, every $\alpha\in Aut(M)$ is unitarily implemented, thus the previous theorem applies and we are done.
\begin{flushright}
$\square$
\end{flushright}

\section{The Hardy algebra compressed to $\mathfrak{M}$.}

In the previous chapter in formula $(\ref{DefnOfGenerSymmetricPart})$ we defined the generalized symmetric $\rho_\pi(H^{\infty}(E))$-coinvariant subspace $\mathfrak{M}$. Remember that we write $\rho_0$ for the representation $\rho_\pi$.
In this chapter we shall prove the reflexivity of the compression of $\rho_0(H^{\infty}(E))$ to the subspace $\mathfrak{M}$, i.e. we shall prove the reflexivity of $Q\rho_0(H^{\infty}(E))|_{\mathfrak M}$, where we write $Q=P_{\mathfrak{M}}$ for the projection onto $\mathfrak{M}$. It should be noted that even in the cases when we know that $\rho_0(H^\infty(E))$ is reflexive, this result does not follows immediately, since the algebras $Alg\ Lat\ (Q\rho_0(H^\infty(E))|_{\mathfrak M})$ and $Q(Alg\ Lat\ \rho_0(H^\infty(E)))|_{\mathfrak M}$ are not the same in general. Note also that the theorem on the reflexivity of $Q\rho_0(H^\infty(E))|_{\mathfrak M}$ generalizes to our setting a result of G. Popescu from ~\cite[Theorem  4.5]{Po6}. I thank Orr Shalit who pointed out that this theorem is also related to the fact that multiplier algebras on a reproducing kernel Hilbert space are always reflexive.

Write $Q_k$ for the projection onto $H(k)=\bigvee_{\eta}\eta^{(k)}(H)$.
Clearly, $Q_k\leq P_k\otimes I_H$. In particular, since $M\otimes_{\pi}H\cong
H=H(0)$, and $E\otimes_{\pi}H=\bigvee_{\eta}\eta(H)=H(1)$, we have
$Q_0=P_0\otimes I_H$ and $Q_1=P_1\otimes I_H$.

Set
$$\tilde{W}_t=\sum_{k\geq 0}e^{ikt}Q_k,$$
and
$$\tilde{\gamma}_t=Ad\tilde{W}_t,$$
where $Ad\tilde{W}_t(T)=\tilde{W}_tT\tilde{W}_t^*$ for $T\in B(\mathfrak{M})$.
Then $\{\tilde{\gamma}_t\}_{t\in \mathbb{R}}$ is an ultraweakly continuous action of
$\mathbb{R}$ on $B(\mathfrak{M})$, which is the gauge automorphism group. As in the previous chapter, the group $\{\tilde{\gamma}_t\}_{t\in \mathbb{R}}$ leaves invariant $Q\rho_0(H^{\infty}(E))|_{\mathfrak{M}}$. The $j$-th Fourier coefficient of $T\in
B(\mathfrak{M})$, associated with $\{\tilde{\gamma}_t\}_{\mathbb{R}}$, is defined by
$$\Phi_j(T)=\frac{1}{2\pi}\int_{0}^{2\pi}e^{-ijt}\tilde{\gamma}_t(T)dt.$$
As for $\mathcal{L}(\mathcal{F}(E))$ this integral ultraweakly converges in $B(\mathfrak{M})$ and leaves invariant
$Q\rho_0(H^{\infty}(E))|_{\mathfrak{M}}$ (since $Q\rho_0(H^{\infty}(E))|_{\mathfrak{M}}$ is $\{\tilde{\gamma}_t\}_{t\in\mathbb{R}}$-invariant).
Let $Y$ be an element of $Alg\ Lat\ Q\rho_0(H^{\infty}(E))|_{\mathfrak M}$. Direct calculation
gives the formula
$$\Phi_j(Y)=\sum_{k\geq 0}Q_{k+j}YQ_k.$$

Here we shall again use the upper index to indicate the degree of the general element of $E^{\otimes s}$, i.e. $\xi^{(s)}\in E^{\otimes s}$.
\begin{lem}
If $Y\in Alg\ Lat\ Q\rho_0(H^{\infty}(E))|_{\mathfrak M}$ then for every
$j$, the operator
$\Phi_j(Y)$ is in $Alg\ Lat\ Q\rho_0(H^{\infty}(E))|_{\mathfrak M}$
\end{lem}
\noindent{Proof.} From $\Phi_j(Y)=\sum_{k\geq 0}Q_{k+j}YQ_k$
we conclude that $\Phi_j(Y)=0$ for $j<0$. Let $j\geq 0$. Note that
$\tilde{W}_t$ has a closed range and by simple calculation we obtain
that $\tilde{W}_t^*=\tilde{W}_{-t}$ and
$$\tilde{W}_t(T_{\xi}\otimes I_H)_{|_{\mathfrak{M}}}=e^{it}Q(T_{\xi}\otimes I_H)\tilde{W}_t|_{\mathfrak{M}}.$$
Indeed, let $x=h_0+\theta_1\otimes h_1+\theta^{(2)}_2\otimes h_2+...+\theta^{(s)}_s\otimes h_s \in \mathfrak{M}$.
Then
$$Q(T_{\xi}\otimes I_H)\tilde{W}_tx=Q(T_{\xi}\otimes I_H)(h_0+e^{it}\theta_1\otimes h_1+
e^{2it}\theta^{(2)}_2\otimes h_2+...+e^{sit}\theta^{(s)}_s\otimes h_s )=$$
$$Q_1(\xi\otimes h_0)+e^{it}Q_2(\xi\otimes\theta_1\otimes h_1)+
e^{2it}Q_3(\xi\otimes\theta^{(2)}_2\otimes h_2)+...+e^{sit}Q_{s+1}(\xi\otimes\theta^{(s)}_s\otimes h_s)=$$
$$=e^{-it}\tilde{W}_t(T_{\xi}\otimes I_H)(h_0+\theta_1\otimes h_1+\theta^{(2)}_2\otimes h_2+...+\theta^{(s)}_s\otimes h_s).$$

It follows that if $\mathcal{M}\in Lat\ Q\rho_0(H^{\infty}(E))|_{\mathfrak M}$ then also
$\mathcal{M}_t:=\tilde{W}_t(\mathcal{M})\in Lat\
Q\rho_0(H^{\infty}(E))|_{\mathfrak M}$.

Let $Y\in Alg\ Lat\ Q\rho_0(H^{\infty}(E))|_{\mathfrak M}$. Then
$$\tilde{\gamma}_t(Y)(\mathcal{M})=\tilde{W}_tY\tilde{W}_{-t}(\mathcal{M})=
\tilde{W}_tY(\mathcal{M}_{-t})
\subseteq
\tilde{W}_t(\mathcal{M}_{-t})=\tilde{W}_t\tilde{W}_{-t}(\mathcal{M})=\mathcal{M}.$$
Integrating, we get $\Phi_j(Y)(\mathcal{M})\subseteq \mathcal{M}$, for all $j$.
\begin{flushright}
$\square$
\end{flushright}
\begin{prop}\label{Restrictions of compressed Phi_j on H} Let $Y\in Alg\ Lat\
Q\rho_0(H^{\infty}(E))|_{\mathfrak M}$ and let us consider the restriction
$\Phi_j(Y)|_H$. Then

1) For $j=0$ there exists $a_0\in M$ such that
$$\Phi_0(Y)|_H=\pi(a_0).$$

2) For $j\geq 1$ there exist $\xi^{(j)}_0\in E^{\otimes j}$ such that
$\Phi_j(Y)|_H=Q_jL^{(0)}_{\xi^{(j)}_0}=Q_j(T_{\xi^{(j)}_0}\otimes I_H)|_H$.
\end{prop}

\noindent{Proof.}
1)  Since $\Phi_0(Y)=\sum_kQ_kYQ_k$, we obtain that $\Phi_0(Y)|_H=Q_0YQ_0$.
Pick some $h\in H$ and set
$$\mathcal{M}_h:=Q\overline{\rho_0(H^{\infty}(E))h}=Q[\overline{\pi(M)h}\oplus
\overline{E\otimes_{\pi}h}\oplus...].$$
Then $\mathcal{M}_h$ is $Q\rho_0(H^{\infty}(E))Q$-invariant, $h\in
\mathcal{M}_h$, and $\Phi_0(Y)h\in Q_0\mathcal{M}_h=\overline{\pi(M)h}$. There
exist a net $\{a_{\iota}\}\subset M$ such that
$\Phi_0(Y)h=\lim_{\iota}\pi(a_{\iota})h$. Let $p\in \pi(M)'$ be any projection.
Then
there are nets $\{b_{\iota}\}$ and $\{c_{\iota}\}$ in $M$ such that
$$\Phi_0(Y)ph=\lim_{\iota}\pi(b_{\iota})ph\,\,\,\,\text{and}\,\,\,\,
\Phi_0(Y)p^{\perp}h=\lim_{\iota}\pi(c_{\iota})p^{\perp}h.$$
Hence,
$$\Phi_0(T)h=\lim_{\iota}\pi(b_{\iota})ph+\lim_{\iota}\pi(c_{\iota})p^{\perp}h,$$
and applying $p$ on both sides, we obtain
$$p\Phi_0(Y)h=\lim_{\iota}\pi(b_{\iota})ph=\Phi_0(Y)ph.$$
It follows that $\Phi_0(Y)$ commutes with $p$ at $h$. Since $h$ is arbitrary we
obtain that $\Phi_0(Y)p=p\Phi_0(Y)$ on $H$. Since this holds for every projection
from the von Neumann algebra $\pi(M)'$, we obtain that $\Phi_0(Y)|_H\in
\pi(M)''=\pi(M)$. Thus,
$\Phi_0(Y)|_{H}=\pi(a_0)$, for some $a_0\in M$.


2) For $j\geq 1$ we have $\Phi_j(Y)(H)=Q_jYQ_0(H)\subset E^{\otimes
j}\otimes_{\pi}H$. As for $\Phi_0(Y)$ we pick an arbitrary $h\in H$ and form the
$Q\rho_0(H^{\infty}(E))Q$-invariant subspace $\mathcal{M}_h=Q\overline{\rho_0(H^{\infty}(E))h}$ as above.
Then $\Phi_j(Y)h\in \mathcal{M}_h$, and $\Phi_j(Y)h\in Q_j\overline{E^{\otimes
j}\otimes h}\subseteq \overline{E^{\otimes j}\otimes h}$. For any projection
$p\in \pi(M)'$ there are
nets $\{\theta^{(j)}_{\iota}\}$, $\{\zeta^{(j)}_{\iota}\}$, $\{\vartheta^{(j)}_{\iota}\}$ in $E^{\otimes j}$
such that
$\Phi_j(Y)h=\lim_{\iota}\theta^{(j)}_{\iota}\otimes h$,
$\Phi_j(Y)ph=\lim_{\iota}\zeta^{(j)}_{\iota}\otimes ph$ and
$\Phi_j(Y)p^{\perp}h=\lim_{\iota}\vartheta^{(j)}_{\iota}\otimes p^{\perp}h$.
Hence
$$\Phi_j(Y)h=\lim_{\iota}\zeta^{(j)}_{\iota}\otimes
ph+\lim_{\iota}\vartheta^{(j)}_{\iota}\otimes p^{\perp}h,$$
and by applying $I_j\otimes p$ on both sides we obtain
$$(I_j\otimes p)\Phi_j(Y)h=(I_j\otimes p)\lim_{\iota}\zeta^{(j)}_{\iota}\otimes
ph=\lim_{\iota}\zeta^{(j)}\otimes ph=\Phi_j(Y)ph.$$
Thus, $$(I_j\otimes p)\Phi_j(Y)=\Phi_j(Y)p$$
on $H$. Since $p$ is an arbitrary projection in $\pi(M)'$ we obtain that
$\Phi_j(Y)$ is intertwines the actions $\iota(\cdot)$ and $I_j\otimes \iota(\cdot)$ of $\pi(M)'$ on $H$ and on $E^{\otimes j}\otimes_\pi H$ respectively, where $\iota$ is the identity action of $\pi(M)'$. By combining with the Fourier
transform $U=U_\pi$ we obtain the operator
$$U\Phi_j(Y):H\rightarrow (E^{\pi})^{\otimes j}\otimes_{\iota}H.$$
Since $U$ also intertwining the actions of $\pi(M)'$ on $E^{\otimes j}\otimes_\pi H$ and on $(E^{\pi})^{\otimes j}\otimes_\iota H$, as it pointed in the Chapter 2, formula $(\ref{Intertwining M' and U})$. Hence there exist the unique
$\xi^{(j)}_0\in E^{\otimes j}$ such that
$$\hat{\xi}^{(j)^*}_0(\eta_1\otimes...\otimes \eta_k\otimes
h)=L^*_{\xi^{(j)}_0}((I_{j-1}\otimes \eta_1)...\eta_k(h)),$$
and in particular,
$$\hat{\xi}^{(j)*}_0(\eta^{\otimes j}\otimes h)=L^*_{\xi^{(j)}_0}((I_{j-1}\otimes
\eta)...\eta(h)).$$
Thus, from the formula $\ref{Formula for hat(xi)}$,
$$U^*\hat{\xi}^{(j)*}_0=L^*_{\xi^{(j)}_0}.$$
Recall also that $\mathfrak{M}=U^*\tilde{\mathfrak{M}}$, and, in particular,
$H(s)=U^*_s\tilde{H}(s)$, $s\geq 0$(in notations of Chapter 4).
Hence, $\Phi_j(Y)|_H=Q_j(T_{\xi^{(j)}_0}\otimes I_H)|_H$.
\begin{flushright}
$\square$
\end{flushright}

\begin{lem} The map $X\mapsto QXQ$ defines a homomorphism from the algebra
$Alg\ Lat\ \rho_0(H^{\infty}(E))$ to the algebra $QAlg\ Lat\
\rho_0(H^{\infty}(E))|_{\mathfrak{M}}$.
\end{lem}
\noindent{Proof.} Enough to show that this map is multiplicative, that is
$QXYQ=QXQYQ$, for every $X,Y\in Alg\ Lat\ \rho_0(H^{\infty}(E))$. But the
coinvariant projection $Q$ can be represented as the difference $I-Q^{\perp}$ of
two invariant projections. Thus, $Q$ is semiinvariant and the map $X\mapsto QXQ$
is multiplicative.
\begin{flushright}
$\square$
\end{flushright}

Let $\eta\in \mathbb{D}(E^{\pi})$. Then $\mathfrak{M}_\eta$ is $Q\rho_0(H^{\infty}(E))Q$-invariant.
To see this, note, as it follows from the preceding lemma, that for every $\xi\in E$ and $a\in M$, $Q(T_\xi\otimes I_H)Q$ and $Q(\phi_{\infty}(a)\otimes I_H)Q$ are the generators of $Q\rho_0(H^{\infty}(E))Q$. Since $Q\mathfrak{M}_\eta=\mathfrak{M}_\eta$ and $\mathfrak{M}_\eta$ is $\rho_0(H^{\infty}(E))$-invariant we deduce that $\mathfrak{M}_\eta$ is $Q\rho_0(H^{\infty}(E))Q$-invariant.

Thus, for every $Y\in Alg\ Lat\ Q\rho_0(H^{\infty}(E))|_{\mathfrak{M}}$ we have $Y^*\mathfrak{M}_{\eta}\subseteq \mathfrak{M}_{\eta}$, hence $\Phi_j(Y)^*\mathfrak{M}_{\eta}\subseteq \mathfrak{M}_{\eta}$ for every $j\geq 0$.

Consider the zeroth coefficient $\Phi_0(Y)$ of a given $Y\in Alg\ Lat\ Q\rho_0(H^{\infty}(E))|_{\mathfrak{M}}$. Take any $h\in H$, then
$$\Phi_0(Y)^*(h+\eta(h)+...+\eta^{(s)}(h)+...)=(k+\eta(k)+...+\eta^{(s)}(k)+...).$$
We have $\Phi_0(Y)^*h=\pi(a_0^*)h=k$ and
$$\Phi_0(Y)^*\eta^{(s)}(h)=\eta^{(s)}(k)=\eta^{(s)}(\pi(a_0^*)h)=
(\phi_{\infty}(a_0^*)\otimes I_H)\eta^{(s)}(h).$$
Hence,
$$\Phi_0(Y)^*(\sum_s\eta^{(s)}(h))=(\phi_{\infty}(a_0^*)\otimes
I_H)(\sum_s\eta^{(s)}(h)),$$
and we obtain that $\Phi_0(Y)^*\mathfrak{M}_{\eta}= (\phi_{\infty}(a_0^*)\otimes
I_H)\mathfrak{M}_{\eta}.$

\begin{thm}\label{Reflexivity of CompressedAlg} The algebra $Q\rho_0(H^{\infty}(E))|_{\mathfrak{M}}$ is reflexive.
\end{thm}
\noindent{Proof.} Every $Y\in Alg\ Lat\
Q\rho_0(H^{\infty}(E))|_{\mathfrak{M}}$ is a ultraweak limit of Cesaro sums of its Fourier coefficients. Hence, enough to show that if $Y\in Alg\ Lat\
Q\rho_0(H^{\infty}(E))|_{\mathfrak{M}}$ then

1)$$\Phi_0(Y)|_{\mathfrak{M}}=(\phi_{\infty}(a_0)\otimes I_H)|_{\mathfrak{M}},$$
and

2)$$\Phi_j(Y)|_{\mathfrak{M}}=Q(T_{\xi^{(j)}_0}\otimes I_H)|_{\mathfrak{M}},$$
for $j\geq 1$.

For 1) let $j=0$. From the equality $\Phi_0(Y)^*\mathfrak{M}_{\eta}=
(\phi_{\infty}(a_0^*)\otimes I_H)\mathfrak{M}_{\eta}$ we conclude that
$$\Phi_0(Y)^*\eta^{(s)}(h)=(\phi_{s}(a_0^*)\otimes I_H)\eta^{(s)}(h),$$
for every $\eta\in E^{\pi}$, $h\in H$ and $s\geq 0$.
Hence,
$$\Phi_0(Y)^*\mathfrak{M}=(\phi_{\infty}(a_0^*)\otimes I_H)\mathfrak{M},$$
and as a consequence
$$\Phi_0(Y)|_{\mathfrak{M}}=(\phi_{\infty}(a_0)\otimes I_H)|_{\mathfrak{M}}.$$
For 2) let
$j\geq 1$. We have $\Phi_j(Y)^*\mathfrak{M}_{\eta}\subseteq
\mathfrak{M}_{\eta}$.
Then $$\Phi_j(Y)^*(\sum_s\eta^{(s)}(h))=\sum_s\eta^{(s)}(k).$$
Clearly, $\Phi_j(Y)^*(\eta^{(s)}(h))=0$ for $0\leq s\leq j-1$.

Set $k=\Phi_j(Y)^*(\eta^{(j)}(h))=(Q_jYQ_0)^*(\eta^{(j)}(h))=
(Q_0Y^*Q_j)(\eta^{(j)}(h))=(T_{\xi^{(j)}_0}^*\otimes I_H)\eta^{(j)}(h)$.
Then, for some $k\in H$,
$$\eta^{(s)}(k)=\eta^{(s)}((T_{\xi^{(j)}_0}^*\otimes
I_H)(\eta^{(j)}(h)))=(T_{\xi^{(j)}_0}^*\otimes I_H)\eta^{(s+j)}(h).$$
Hence,
$$\Phi_j(Y)^*|_{\mathfrak{M}_{\eta}}=(T_{\xi^{(j)}_0}^*\otimes
I_H)|_{\mathfrak{M}_{\eta}},$$
and as in 1),
$$\Phi_j(Y)|_{\mathfrak{M}}=Q(T_{\xi^{(j)}_0}\otimes I_H)|_{\mathfrak{M}}.$$
Thus, every $\Phi_j(Y)$, $j\geq 0$ is in the compressed algebra $Q\rho_0(H^{\infty}(E))|_{\mathfrak{M}}$, which implies that this algebra is reflexive.
\begin{flushright}
$\square$
\end{flushright}

~\\
\textbf{Acknowledgement.} I thank Baruch Solel for suggesting problems, constant help and encouragement.

\end{document}